\documentclass[a4paper,12pt,onecolumn]{article}

\usepackage{hyperref}
\usepackage[vmargin=2cm,hmargin=2cm,headheight=14.5pt,top=2cm,headsep=.5cm]{geometry}

\usepackage{bm}
\usepackage[utf8]{inputenc}

\usepackage[autobold]{mathfixs}
\usepackage{empheq}
\usepackage{stackrel}
\usepackage{cases}
\usepackage{float}
\usepackage{mathtools}
\usepackage{amsthm,amsmath,amscd}
\usepackage{makeidx}
\usepackage[charter]{mathdesign}
\usepackage{perpage}
\usepackage{bm}
\usepackage{tikz-cd}
\tikzcdset{every label/.append style = {font = \small}}
\usepackage{caption}
\usepackage{xcolor}
\usepackage{graphicx}
\DeclareMathSizes{12}{12}{8}{6}

\usepackage{cite}
\usepackage{url}

\usepackage[ddmmyy]{datetime}
\usepackage{enumitem}
\setlist[itemize,2]{label=$\centerdot$}
\setlist[itemize,3]{label=$\triangle$}
\usepackage{framed}
\usepackage[symbol]{footmisc}
\usepackage{marvosym}
\usepackage{perpage}
\MakePerPage[1]{footnote}

\makeatletter
\renewcommand*{\@makefnmark}{\hbox{\@textsuperscript{%
			\normalfont\@thefnmark}}}
\makeatother

\newtheoremstyle{ptheorem}{1em}{0em}{\itshape}{}{\bfseries}{.}{.5em}{\thmname{#1}\thmnumber{
		#2}\thmnote{ (\hspace{-.01pt}{#3})}}

\theoremstyle{ptheorem}

\newtheorem{thm}{Theorem}[section]
\newtheorem{pro}[thm]{Proposition}
\newtheorem{lem}[thm]{Lemma}
\newtheorem{cor}[thm]{Corollary}

\newtheoremstyle{hdef}{1em}{0em}{}{}{\bfseries}{.}{.5em}{\thmname{#1}\thmnumber{
		#2}\thmnote{ (\hspace{-.01pt}{#3})}}
\theoremstyle{hdef}

\newtheorem{dfn}[thm]{Definition}
\newtheorem{rem}[thm]{Remark}

\numberwithin{equation}{section}
\numberwithin{figure}{section}

\renewcommand{\k}{\kappa}
\renewcommand{\phi}{\varphi}

\newcommand{\olb}[1]{%
	\vbox{\offinterlineskip\ialign{\hfil##\hfil\cr
			$\rotatebox[origin=c]{90}{$]$}$\cr\noalign{\kern-.45ex}{$#1$}\cr}}}

\newcommand{\noop}[1]{}

\usepackage{stmaryrd}

\allowdisplaybreaks

\begin{document}

	\title{Compactness criteria for Stieltjes function spaces and applications}

\date{}

\author{
	F. Javier Fernández$^*$\\
	\small e-mail: fjavier.fernandez@usc.es\\
	F. Adri\'an F. Tojo$^*$ \\
	\small e-mail: fernandoadrian.fernandez@usc.es
\\Carlos Villanueva\\ \small e-mail: carlos.villanuevamariz@lmh.ox.ac.uk}

\maketitle

\begin{abstract}
In this work we study some topological aspects of function spaces arising in Stieltjes differential calculus. Chief among them are compactness results related to the Ascoli-Arzelà and Kolmogorov-Riesz theorems, as well as their applications to Stieltjes-Sobolev spaces and decomposable functions.
\end{abstract}



\footnotetext{\noindent$^*$\emph{CITMAga, 15782, Santiago de Compostela, Spain},\\ \emph{Departamento de Estatística, Análise Matemática e Optimización},\\ \emph{Universidade de Santiago de Compostela, 15782, Facultade de Matemáticas, Santiago, Spain.}}

\section{Introduction}
Stieltjes differential calculus is a powerful tool in the study of complex differential problems \cite{LoRo14}. It is based on a notion of derivative with respect to a left continuous and non-decreasing function that we call a derivator --see Definition~\ref{g_der}. This theory has been proven particularly useful while modeling processes in which the features studied have sudden changes in time --see, for instance \cite{FriLo17} and \cite{MarquezTesis}.\par

On the other hand, compactness results are crucial when studying the existence of solutions of differential problems (for example, if we want to use fixed point index methods). The first of this kind of results for the case of Stieltjes function spaces appeared in \cite[Theorem~4.2]{FriLo17}, where a sufficient condition of relative compactness in $\operatorname{C}_g([a,b],{\mathbb R})$ was given. This same result was deepened (in a slightly more general setting) in \cite{FriTo20}, giving a characterization for compactness in $\operatorname{BC}_g([a,b],{\mathbb R})$ (Theorem~\ref{thmbc}). Nonetheless, the theory is still young, with many classical results without a counterpart in this setting and various open questions unanswered; among them, for instance, if the --rather unwieldy-- condition of $g$-stability --Definition~\ref{dgs}-- can be restated in a more meaningful way, and whether it can be related to the uniform $g$-equicontinuity.\par
The aim of this work is, therefore, to shed light on these issues, answering some of the open questions, providing new compactness results which are applicable to the study of the existence of solution of Stieltjes differential problems and offering some applications of the theory developed.

This work is structured as follows. First, in Section 2, we provide a brief summary of preliminary concepts for the convenience of the reader. Then, we show some interesting topological aspects of the space of $g$-continuous functions in Section 3, where we also relate the results presented to those in the literature. Furthermore, we apply the new results to obtain a version of the Weierstrass approximation theorem for uniformly $g$-continuous functions.\par
In Section 4, we move on to spaces of $g$-integrable functions. There, some results from \cite{Fernandez2022} are combined with our new findings in order to provide an extension for Stieltjes measures of the Kolmogorov-Riesz compactness theorem.\par
In Sections 5 and 6 we showcase different applications of the theory. On the one hand, we derive some compactenss results for Stieltjes-Sobolev spaces and, on the other, we provide classification and compactenss results for decomposable functions.

\section{Preliminaries}
Throughout the article, ${\mathbb F}$ will be ${\mathbb R}$ or ${\mathbb C}$ and $g:\mathbb{R}\to\mathbb{R}$ will be a non-decreasing and left \mbox{continuous} function, which we will refer to as \emph{derivator}. For a given derivator, we define the sets:
\[ \begin{aligned}
&C_g:=\{t\in\mathbb{R}:g\text{ is constant on } (t-\varepsilon,t+\varepsilon)\text{ for some }\varepsilon>0\},\\
&D_g:=\{t\in\mathbb{R}:\Delta g(t)>0\},
\end{aligned}\] 
where $\Delta g(t)=g(t^+)-g(t)$, and $g(t^+)$ is the right hand side limit of $g$ at $t$. $C_g$ is open in the usual topology, so it can be expressed as 
\begin{equation}\label{c_g}
C_g=\bigcup_{n\in\Lambda}(a_n,b_n),
\end{equation}
where $\Lambda$ is countable and the union is disjoint. With this in mind, we also consider the sets
\[ N_g^-:=\{a_n:n\in\Lambda\}\backslash D_g,\ N_g^+:=\{b_n:n\in\Lambda\}\backslash D_g,\ N_g=N_g^-\cup N_g^+.\] 
By $\mu_g$ we will denote the Stieltjes measure associated to $g$ (see, for example, \cite{Athreya2006} for details). A set or a function is $g$-\emph{measurable} if it is measurable with respect to the complete $\sigma$-algebra obtained through Carathéodory's extention theorem \cite[Theorem~1.3.3]{Athreya2006}. We denote this $\sigma$-algebra by $\mathcal{M}_g$. ${L}^1_g(X;\mathbb{F})$ is the set of equivalence classes of $\mu_g$-integrable functions on a $g$-measurable set $X$ with values in $\mathbb{F}$ where $f\equiv h$ if and only if $\mu_g((f-h)^{-1}({\mathbb F}\backslash\{0\}))=0$, i.e., $f=h$ $\mu_g$-a.e. Given $f\in {L}^1_g(X;\mathbb{F})$, we denote its integral by
\[ \int_X f(s) \operatorname{d} \mu_g(s),\ f\in{L}^1_g(X;\mathbb{F}).\] 
As usual, we do not make the difference between the equivalence classes in ${L}^1_g(X;\mathbb{F})$ and their representatives.
\begin{dfn}
	A function $f:[a,b]\to\mathbb{R}$ is $g$-\emph{continuous} at a point $t\in[a,b]$ if for every $\varepsilon>0$ there exists $\delta>0$ such that 
	\[ \left\lvert f(t)-f(s)\right\rvert<\varepsilon,\ \text{ for all }s\in[a,b],\ \left\lvert g(t)-g(s)\right\rvert<\delta. \] 
	$g$-continuity on $[a,b]$ and uniform $g$-continuity are defined analogously to the usual case. We denote by $\operatorname{C}_g([a,b];\mathbb{F})$ the \emph{set of} $g$-\emph{continuous functions} on $[a,b]$, and $\operatorname{BC}_g([a,b],\mathbb{F})$ the \emph{set of} \emph{bounded} $g$-\emph{continuous functions} on $[a,b]$.
\end{dfn}	
\begin{rem}
	 For an arbitrary $g$, it is not always the case that $\operatorname{BC}_g([a,b],\mathbb{F})=\operatorname{C}_g([a,b],\mathbb{F})$. That is, $g$-continuous functions on compact intervals are not necessarily bounded (see \cite[Example~3.19]{MarquezTesis}). 
\end{rem}
The following result describes properties satisfied by $g$-continuous functions.
\begin{pro}[{\cite[Proposition 3.2]{FriLo17}}]\label{profg}
	If $f:[a,b]\to\mathbb{F}$ is $g$-continuous on $[a,b]$, then:
	\begin{itemize}
		\item $f$ is continuous from the left at every $t\in(a,b]$;
		\item if $g$ is continuous at $t\in[a,b)$, then so is $f$;
		\item if $g$ is constant on some $[\alpha,\beta]\subset[a,b]$, then so is $f$.
	\end{itemize}
\end{pro}
The concept of $g$-continuity is closely related to $g$-topologies, which we now define.
\begin{dfn}
For a derivator $g:\mathbb{R}\to\mathbb{R}$, the topology induced on the real line by the pseudometric
\[ \rho(x,y)=\left\lvert g(x)-g(y)\right\rvert\] 
is called $g$-\emph{topology}, and we denote it by $\tau_g$.
\end{dfn}
\begin{rem}
	Analogously to the standard case, it can be seen that $g$-continuity of $f:\mathbb{R}\to\mathbb{F}$ is equivalent to the continuity of the map between topological spaces $f:(\mathbb{R},\tau_g)\to(\mathbb{F},\tau_{usual})$.
\end{rem}
Next we define the central object of Siteltjes differential calculus.
\begin{dfn}\label{g_der}
We define the \emph{Stieltjes derivative}, or $g$-\emph{derivative}, of $f:[a,b]\to\mathbb{F}$ at a point $t\in[a,b]$ as
\[ f_g^{\prime}(t)=\left\{\begin{aligned}
	&\lim _{s \rightarrow t} \frac{f(s)-f(t)}{g(s)-g(t)}, && t \notin D_g \cup C_g, \\
	&\lim _{s \rightarrow t^{+}} \frac{f(s)-f(t)}{g(s)-g(t)}, && t \in D_g, \\
	&\lim _{s \rightarrow b_n^{+}} \frac{f(s)-f\left(b_n\right)}{g(s)-g\left(b_n\right)}, && t \in C_g, t \in\left(a_n, b_n\right),
\end{aligned}\right.\] 
where $a_n$ and $b_n$ are as in~\eqref{c_g}, provided the corresponding limit exists, in which case we say that $f$ is $g$-\emph{differentiable} at $t$.
\end{dfn}

We introduce now the concept of absolutely $g$-continuous function, along with a version of the fundamental theorem of calculus for the $g$-derivative. This is one of the pillars of the theory.
\begin{thm}[{\cite[Theorem 5.4]{LoRo14}}]\label{tfc}
	Let $F:[a,b]\to\mathbb{F}$. The following are equivalent:
	\begin{enumerate}
		\item F is \emph{abolutely} $g$-\emph{continuous} on $[a,b]$, that is, for every $\varepsilon>0$ there exists $\delta>0$ such that for any family of pairwise disjoint subintervals $\{(c_n,d_n)\}_{n=1}^{m}$,
		\[ \sum_{n=1}^m\left(g\left(b_n\right)-g\left(a_n\right)\right)<\delta \Longrightarrow \sum_{n=1}^m\left\lvert F\left(b_n\right)-F\left(a_n\right)\right\rvert<\varepsilon .\]  
		\item $F$ satisfies:
		\begin{enumerate}
			\item $F_g'(t)$ exists for $\mu_g$-a.e. $t\in[a,b)$;
			\item $F_g'\in{L}_g^1([a,b),\mathbb{F})$;
			\item for each $t\in[a,b]$,
			\[ F(t)=F(0)+\int_{[0,t)}F_g'(s)\mathrm{d}\mu_g(s).\] 
		\end{enumerate}
	\end{enumerate}

\end{thm}	

Lastly, we review some facts about compactness in metric spaces.
\begin{dfn}
	Let $(X,\tau)$ be a topological space. We say that $U\subset X$ is \emph{relatively compact} if $\overline{U}$ (the closure of $U$) is compact.
\end{dfn}
\begin{dfn}
	Let $(M,d)$ be a metric space. $U\subset M$ is \emph{totally bounded} if for every $\varepsilon>0$ there exists a finite set $F_{\varepsilon}\subset U$ such that $U\subset \bigcup_{x\in F_{\varepsilon}}B_M(x,\varepsilon)$. $F_\varepsilon$ is called an $\varepsilon$-\emph{net}.
\end{dfn}
Instead of relative compactness, we will often check total boundedness. The following proposition establishes that these are equivalent for complete metric spaces.
\begin{pro}[{\cite[Theorems 7.5 and 8.2, Exercise 15 p. 110]{Carothers2000}}]\label{equiv_rel_compacto}
	Let $(M,d)$ be a metric space, $A\subset M$. The following are equivalent:
	\begin{enumerate}
		\item $A$ is totally bounded;
		\item Every sequence in $A$ has a Cauchy subsequence.
	\end{enumerate}
	Furthermore, if $(M,d)$ is complete, the previous conditions are equivalent to
	\begin{enumerate}[resume]
		\item $A$ is relatively compact.
	\end{enumerate} 
\end{pro}

\section{Compactness in the space $\operatorname{BC}_g([a,b],{\mathbb F})$}
An incorrect proof of the second countability of $g$-topologies appears in \cite{FriLo17}. The mistake was pointed out in \cite{MarquezTesis}, but without giving an alternative proof. That is the first result of the section.
\begin{pro}\label{2do_num}
	The set
	\begin{equation}\label{conj_num}
		\mathcal{C}=(\{(a,b):a\in\mathbb{Q}\cup D_g\cup N_g^{+}\text{ y } b\in\mathbb{Q}\cup D_g\cup N_g^{-}\}\cup\{(c,d]:c\in\mathbb{Q}\cup D_g\cup N_g^{+},d\in D_g\})\cap\tau_g
	\end{equation}
	is a base of $\tau_g$. $g$-topologies are therefore second countable, and also $\tau_g\subset\mathcal{B}(\mathbb{R})$ (the Borel $\sigma$-algebra of $\mathbb{R}$). 
\end{pro}
\begin{proof}
	We start by making sure that the sets of the form $(c,d]\in\tau_g$ (with $d\in D_g$ and $c$ any real number) are the countable union of sets in $\mathcal{C}$. If $c\in D_g$, we are done and, if $c\notin D_{g}$, we necessarily have that $(c,c+1/n)\not\subset C_g$ for any $n$ (otherwise $(c,d]\notin\tau_g$). Taking this into account, there exists a sequence $\{c_n\}_{n\in\mathbb{N}}\subset(c,d]\cap(\mathbb{Q}\cup D_g\cup N_g^{+})$ converging to $c$, and such that $f(c_n)<f(t)$ for every $t>c_n$. Under this conditions, $(c_n,d]\in\mathcal{C}$, y $(c,d]=\bigcup_{n\in\mathbb{N}}(c_n,d]$. A similar argument is valid for $(a,b)\in\tau_g$, with $a,b\in\mathbb{R}$.\par
	Any element of $\tau_g$ is the union (possibly arbitrary) of open $g$-balls (that is, sets of the form $\{t\in [a,b]\ :\ \left\lvert g(t)-g(t_0)\right\rvert<r\}$ where $t_0\in[a,b]$ and $r>0$), which are always of the form
	\[ g^{-1}((\alpha, \beta))=\{t \in \mathbb{R}: \alpha<g(t)<\beta\}, \quad \text { with } \alpha, \beta \in \mathbb{R}, \alpha<\beta.\] 
	It is easy to see that $g^{-1}((\alpha,\beta))=(a,b]$ if $b\in D_g$ and $g(b)<\beta\leqslant g(b^{+})$, and that $g^{-1}((\alpha,\beta))=(a,b)$ otherwise. 	Therefore, if $U\in\tau_g$, there exist two (arbitrary) families of open balls of $\tau_g$, $\{(a_i,b_i)\}_{i\in\mathcal{I}}$ and $\{(a_j,b_j]\}_{j\in\mathcal{J}}$ with $b_j\in D_g$, such that
	\[ U=\bigcup_{i\in\mathcal{I}}(a_i,b_i)\cup\bigcup_{j\in\mathcal{J}}(a_i,b_i].\] 
	With either Sorgenfrey's topology \cite[p. 79]{Hart2004} or the usual one, $\mathbb{R}$ is hereditarily Lindelöf, so both unions admit a countable subcover. That is,
	\[ U=\bigcup_{i\in I}(a_i,b_i)\cup\bigcup_{j\in J}(a_i,b_i],\] 
	where $I\subset\mathcal{I}$ and $J\subset\mathcal{J}$ are countable. $U$ is then the countable union of sets that are themselves the countable union (by what was shown at the beginning) of elements in $\mathcal{C}$, and we are done.
\end{proof}

The characterization of compact sets in the space $\operatorname{BC}_g([a,b],{\mathbb F})$ was studied in \cite{FriTo20} (although in \cite{FriTo20} the functions are assumed to take real values, that does not change the validity of the proofs). In order to present it we need the following definitions.
\begin{dfn} A set $S\subset \operatorname{BC}_g([a,b],{\mathbb F})$ is said to be \emph{$g$-equicontinuous} if, for every $\varepsilon\in{\mathbb R}^+$ and $t\in {[a,b]}$, there exists $\delta\in\mathbb{R}$ such that $\left\lvert f(t)-f(s)\right\rvert<\varepsilon$ for every $f\in S$ and every $s\in [a,b]$ such that $\left\lvert g(t)-g(s)\right\rvert<\delta$. We say that $S$ is \emph{uniformly $g$-equicontinuous} if, for every $\varepsilon\in{\mathbb R}^+$, there exists $\delta\in{\mathbb R}^+$ such that $\left\lvert f(t)-f(s)\right\rvert<\varepsilon$ for every $f\in S$ and every $t,s\in [a,b]$ such that $\left\lvert g(t)-g(s)\right\rvert<\delta$.
\end{dfn}

\begin{dfn}\label{dgs}
	A set $S\subset \operatorname{BC}_g([a,b],{\mathbb F})$ is said to be \emph{$g$-stable} if, for every $t\in[a,b)\cap D_g$ and every $\varepsilon\in{\mathbb R}^+$, there exist $\delta\in{\mathbb R}^+$ and a finite covering $\{X_k\}_{k=1}^m$ of $[t,t+\delta)\cap I$ such that $\left\lvert f(x)-f(y)\right\rvert<\varepsilon$ for every $f\in S$, every $x,y\in X_k$ and every $k=1,\dots,m$.
\end{dfn}

With these definitions we get the following result.
\begin{thm}[{\cite[Theorem~4.14]{FriTo20}}]\label{thmbc}
	Let $S\subset \operatorname{BC}_g([a,b],{\mathbb F})$, then $S$ is precompact if and only if
	\begin{enumerate}
		\item $S(t)$ is bounded for all $t\in [a,b]$;
		\item $S$ is $g$-equicontinuous;
		\item $S$ is $g$-stable.
	\end{enumerate}
\end{thm}

Furthermore, we can simplify the hypotheses in Theorem~\ref{thmbc} by observing the following.
\begin{lem}[{\cite[Lemma 4.13]{FriTo20}}] \label{unif-eq=stab}
	Let $S\subset \operatorname{BC}_g([a,b],{\mathbb F})$. If $S$ is uniformly $g$-equicontinuous,
	then $S$ is $g$-equicontinuous and $g$-stable.
\end{lem}
Therefore, we have the following corollary.

\begin{cor}[{\cite[Theorem~4.15]{FriTo20}}]\label{cor:precompact}
	Let $S\subset \operatorname{BC}_g([a,b],{\mathbb F})$. If
	\begin{enumerate}
		\item $S(t)$ is bounded for all $t\in I$ and
		\item $S$ is uniformly $g$-equicontinuous,
	\end{enumerate}
	then $S$ is precompact.
\end{cor}

\cite{FriTo20} leaves open the question of whether this corollary is really an equivalence, that is, if $g$-equicontinuity and $g$-stability imply uniform $g$-equicontinuity, which would make the statement in Lemma~\eqref{unif-eq=stab} an equivalence. We proceed to settle this issue with the following result.

\begin{pro}\label{equi_unif_equi} 
	Let $F\subset{\operatorname{BC}}_g([a,b],\mathbb{F})$ be $g$-equicontinuous. Then the following statements are equivalent:
	\begin{enumerate}
		\item $F$ is uniformly $g$-equicontinuous.
		\item For every $t\in[a,b)\cap D_g$ and $f\in F$, there exists $\lim_{s\to t^{+}}f(s)\in\mathbb{R}$ uniformly on $F$. 
	\end{enumerate}
\end{pro}
\begin{proof} 	 
	$(1)\Rightarrow(2)$ Suppose in the first place that the second statement does not hold. Then, for some $t\in[a,b)\cap D_g$, there exist $\alpha\in\mathbb{R}^{+}$, $\{f_n\}_{n\in\mathbb{N}}\subset F$ and sequences $\{x_n\}_{n=1}^{\infty},\{y_n\}_{n=1}^{\infty}\subset(t,b]$ converging to $t$ such that $\left\lvert f_n(x_n)-f_n(y_n)\right\rvert\geq\alpha>0$. $g$ is increasing, and therefore regulated, so $\left\lvert g(x_n)-g(y_n)\right\rvert\rightarrow 0$ when $n\rightarrow \infty$. Consequently, $F$ can not be uniformly $g$-equicontinuous.\par 
	$(2)\Rightarrow(1)$ Assume that $2$ holds, and fix $\varepsilon\in\mathbb{R}^{+}$. 
	\begin{itemize}
		\item If $t\in[a,b]\backslash D_g$, there exists $\delta_t>0$ with $\sup_{f\in F}\left\lvert f(s)-f(t)\right\rvert<\frac{\varepsilon}{4}$ if $\left\lvert g(s)-g(t)\right\rvert<\delta_t$. Thanks to the continuity of $g$ at $t$, there also exists $\eta(\delta_t)>0$ with $\left\lvert g(s)-g(t)\right\rvert<\delta_t$ for every $s\in(t-\eta(\delta_t),t+\eta(\delta_t))$. Define $I_t:=(t-\eta(\delta_t),t+\eta(\delta_t))\cap[a,b]$.
		\item If $t\in[a,b)\cap D_g$, by hypothesis, there exists $\delta_t^{1}$ with $\sup_{f\in F}\left\lvert f(x)-f(y)\right\rvert<\frac{\varepsilon}{4}$ if $x,y\in(t,t+\delta_t^{1})$. There also exists $\delta_t^{2}>0$ such that $\sup_{f\in F}\left\lvert f(s)-f(t)\right\rvert<\frac{\varepsilon}{4}$ if $\left\lvert g(s)-g(t)\right\rvert<\delta_t^2$, and $\eta(\delta_t^{2})$ with $\left\lvert g(s)-g(t)\right\rvert<\delta_t^2$ for $s\in(t-\eta(\delta_t^{2}),t]$ (thanks to $g$ being left continuous). Define $I_t:=(t-\eta(\delta_t^{2}),t+\delta_t^{1	})\cap[a,b]$.
	\end{itemize}
	$\{I_t\}_{t\in[a,b]}$ is an open cover of $([a,b],\tau_u)$, so we can obtain a finite subcover $\{I_{t_i}\}_{i=1}^{n}\cup\{I_{\widehat{t}_j}\}_{j=1}^{m}$, where $t_i\in[a,b]\backslash D_g$ and $\widehat{t}_j\in[a,b]\cap D_g$. We can assume, discarding redundant intervals, that none of them is contained in the union of the rest, and that they are ordered so that
	\begin{equation}\label{prop_intervalos}
		I_i\cap I_{i+1}\neq\emptyset \text{ and } I_i\cap I_{i+k}=\emptyset \text{ if } k\geq 2.
	\end{equation}\par
	Now for an interval $I$, we will denote by $V(g,I)$ the total variation of $g$ on $I$. Consider
	\[ \delta:=\min\left\{\left\{\frac{V(g,I_s)}{2}\ :\ s\in\{t_1,...,t_n,\widehat{t}_1,...,\widehat{t}_m\},\ V(g,I_s)>0\right\}
	\cup\left\{\frac{g(\widehat{t}_j^{+})-g(\widehat{t}_j)}{2}:j=1,...,m\right\}\right\}.\] 
	We have that $0<\delta<\infty$ as $g$ is non decreasing and, thus, $V(g,[a,b])=g(b)-g(a)$. Let us fix any $f\in F$ and see that it satisfies
	\[ \left\lvert f(s)-f(t)\right\rvert<\varepsilon \text{ whenever } s,t\in[a,b],\ s<t \text{ are such that }\left\lvert g(t)-g(s)\right\rvert=
	V(g,[t,s])<\delta.\] \par
	By the definition of $\delta$, there can be no $\widehat{t}_i\in D_g$ in the interval $[t,s)$. There exist $I_{i_1},...,I_{i_n}$ satisfying~\eqref{prop_intervalos} such that $t\in I_{i_1}$ y $s\in I_{i_n}$; we can also suppose, discarding intervals from the beginning and the end, that $I_{i_2},...,I_{i_{n-1}}\subset(s,t)$ (the remaining cases are simpler, because they involve only one or two intervals).\par
	Since $\left\lvert g(s)-g(t)\right\rvert<\delta$, again by the definition of $\delta$, $V(g,I_{i_2})=\cdots=V(g,I_{i_{n-1}})=0$ and thanks to~\eqref{prop_intervalos}, $g$ is constant on $\bigcup_{k=2}^{n-1}I_{i_k}$. $g$-continuity implies that $f$ is also constant on $I_{i_{j}}$ for every $j=2,...,n-1$ --see Proposition~\ref{profg}, so if we choose $t_{i_j,i_{j+1}}\in I_{i_{j}}\cap I_{i_{j+1}}$ for each $j=1,...,n-1$, we have that 
\begin{align*}
		\left\lvert f(t)-f(s)\right\rvert \leqslant &\left\lvert f(t)-f(t_{i_1})\right\rvert+\left\lvert f(t_{i_1})-f(t_{i_1,i_2})\right\rvert+\left\lvert f(t_{i_1,i_2})-f(t_{i_2})\right\rvert+\cdots\\
		&+\left\lvert f(t_{i_{n-1}})-f(t_{i_{n-1},i_n})\right\rvert+\left\lvert f(t_{i_{n-1},i_n})-f(t_{i_n})\right\rvert+\left\lvert f(t_{i_n})-f(s)\right\rvert\\
		<&\frac{\varepsilon}{4}+\frac{\varepsilon}{4}+0+\cdots+0+\frac{\varepsilon}{4}+\frac{\varepsilon}{4}=\varepsilon.\qedhere
	\end{align*}

\end{proof}

\begin{cor}\label{unif_reglada}
	For a function $f\in\operatorname{BC}_g([a,b],\mathbb{F})$, the following are equivalent
	\begin{enumerate}
		\item $f$ is uniformly $g$-continuous,
		\item $f$ is regulated.
	\end{enumerate}
\end{cor}
\begin{proof}
	It is enough to take $S=\{f\}$ in Proposition~\ref{equi_unif_equi}. 
\end{proof}
\begin{rem} Corollary~\ref{unif_reglada} shows that
	the equivalence in Lemma~\eqref{unif-eq=stab} does not hold. As a counter example it is enough to take
	\[ g(t)= \begin{dcases}t, & t\leqslant 0, \\ t+1, & t > 0,\end{dcases}\qquad f(t)= \begin{dcases}t, & t\leqslant 0, \\ \sin\left(\frac{1}{t}\right), & t > 0,\end{dcases}\] 
	and $S=\{f\}$. $S$ is $g$-equicontinuous and $g$-stable, but not uniformly $g$-equicontinuous. \cite[Example 3.3]{FriLo17} also illustrates this point.
\end{rem}
\begin{rem}
	In Corollary~\ref{unif_reglada}, since $f\in\operatorname{BC}_g([a,b],\mathbb{F})$ is continuous at each point of $[a,b]\backslash D_g$, and left continuous in $[a,b]$, $f$ is regulated if and only if its right-hand limit exists for every $t\in D_g$.
\end{rem}
\begin{rem}\label{rem_ascoli}
	An independent proof of Corollary~\ref{cor:precompact} would not have requiered a version of Ascoli-Arzelà as general as the one involved in Theorem~\ref{thmbc}. Taking Proposition~\ref{equi_unif_equi} into account, a family of functions of ${\operatorname{BC}}_g([a,b],{\mathbb F})$ that satisfies the conditions of Corolary~\ref{cor:precompact} is under the hypotheses of \cite[Theorem 1]{Hildebrandt1966}.
\end{rem}

Corollary~\ref{unif_reglada} has interesting applications. Among them, it offers a way to prove a version of the Weierstrass Approximation Theorem for the space $\operatorname{BUC}_g({\mathbb R},{\mathbb F})$ of bounded uniformly $g$-continuous functions with the supremum norm. Given that the uniform limit of uniformly $g$-continuous functions is a uniformly $g$-continuous function (this is due to Corollary~\ref{unif_reglada} and the fact that the uniform limit of $g$-continuous functions is $g$-continuous \cite[Theorem 3.4]{FriLo17} and the uniform limit of regulated functions is regulated), $\operatorname{BUC}_g({\mathbb R},{\mathbb F})$ is a Banach subspace of the Banach space of bounded $g$-continuous functions with the supremum norm $\operatorname{BC}_g({\mathbb R},{\mathbb F})$.

\begin{thm}\label{thmbc2}
		Let $S\subset \operatorname{BUC}_g([a,b],{\mathbb F})$, then $S$ is precompact if and only if
		\begin{enumerate}
			\item $S(t)$ is bounded for all $t\in [a,b]$;
			\item $S$ is uniformly $g$-equicontinuous.
		\end{enumerate}
	\end{thm}
\begin{proof}
Those two conditions implying precompactness is precisely the statement of Corollary~\ref{cor:precompact} (notice that since $S$ is uniformly $g$-equicontinuous, we have in particular that $S\subset\operatorname{BUC}_g([a,b])$).\par 
For the other direction, Theorem~\ref{thmbc} tells us that $S(t)$ is bounded for each $t\in[a,b]$, and that $S$ is $g$-equicontinuous. Since $S$ is a family of regulated functions precompact in the supremum norm topology, \cite[Theorem 1]{Hildebrandt1966} tells us that the right hand side limits at the points of $D_g$ exist uniformly on $S$. Applying Proposition~\ref{equi_unif_equi}, $S$ is uniformly $g$-equicontinuous.
	\end{proof}
\begin{thm}\label{thmcomp} Let $f\in \operatorname{BUC}_g({\mathbb R},{\mathbb F})$. Then $f=\widetilde f\circ g$ where $\widetilde f$ is a uniformly continuous function.
\end{thm}
\begin{proof}
	Let $c=\sup g({\mathbb R})\in(-\infty,\infty]$. We start by defining the function $\sigma (t)=\sup g^{-1}(t)$ for every $t\in(-\infty,c)\cap g({\mathbb R})$. If $c\in g({\mathbb R})$, we define $\sigma (c)=t$ for some $t\in{\mathbb R}$ such that $g(t)=g(c)$. This way, we have defined $\sigma $ on $g({\mathbb R})$. 	Since $g$ is increasing and left continuous, it is lower semicontinuous, so $\sigma (g(t))=\sup\{s\in{\mathbb R}\ :\ g(s)\leqslant g(t)\}\in g({\mathbb R})$ and $g(\sigma (g(t))=g(t)$ for every $t\in{\mathbb R}$.

	Let $h=f\circ \sigma $ on $g({\mathbb R})$. Then, for $t\in{\mathbb R}$, $h(g(t))=f(\sigma (g(t))$. Since $g(\sigma (g(t))=g(t)$ and $f$ is $g$-continuous, we have that $f(\sigma (g(t)))=f(t)$ --see \cite[Lemma 2.15]{Fernandez2022a}. We conclude that $h\circ g=f$.

	Now, $h$ is continuous. Indeed, given $\varepsilon\in{\mathbb R}^+$, since $f$ is uniformly $g$-continuous, there exists $\delta\in{\mathbb R}^+$ such that, if $\left\lvert g(s)-g(t)\right\rvert<\delta$, then $\left\lvert f(s)-f(t)\right\rvert<\varepsilon$. Given that $h\circ g=f$, we conclude that if $\left\lvert g(s)-g(t)\right\rvert<\delta$ then $\left\lvert h(g(s))-h(g(t))\right\rvert<\varepsilon$, so $h$ is uniformly continuous on $g({\mathbb R})$.

	Furthermore, since $f$ is uniformly $g$-continuous, it is regulated by Corollary~\ref{unif_reglada}. $\sigma $ is increasing, so the composition $f\circ \sigma =h$ is regulated as well. Observe that ${\mathbb R}\backslash g({\mathbb R})$ has no isolated points because $g$ is non-decreasing and, therefore, the connected components of ${\mathbb R}\backslash g({\mathbb R})$ are nondegenerate intervals. In fact, we can write \[ {\mathbb R}\backslash g({\mathbb R})=\bigcup_{d\in D_g}(g(d),g(d^+)]\cup (-\infty,\inf g({\mathbb R})]\cup [\sup g({\mathbb R}),\infty),\] 
	in the case $g$ is bounded, jump points are not followed by constancy intervals and $g$ does not reach a maximum or minimum, with similar expressions (the inclusion of the endpoints of the intervals may vary) depending on other conditions such as whether $g$ is unbounded from above or below and the placement of its intervals of constancy.
	This fact implies that $h$ can be continuously extended, in a unique way, to a uniformly continuous function $\widetilde h$ defined on $\overline{g({\mathbb R})}$, as the lateral limits of $h$ exist at every point of $\overline{g({\mathbb R})}$ where they can be considered due to the fact that it is regulated. Finally, $\overline{g({\mathbb R})}$ is closed, so its complement is open (in the usual topology) and ${\mathbb R}\backslash \overline{g({\mathbb R})}=\bigcup_{n\in\Gamma}(c_n,d_n)$ where $\Gamma\subset{\mathbb N}$, $c_n,d_n\in[-\infty,\infty]$ and the intervals $(c_n,d_n)$ are pairwise disjoint. For instance, in the same situation as above, we would have
	\[ {\mathbb R}\backslash \overline{ g({\mathbb R})}=\bigcup_{d\in D_g}(g(d),g(d^+))\cup (-\infty,\inf g({\mathbb R}))\cup (\sup g({\mathbb R}),\infty).\] 
	We define
	\[ \widetilde f(x):=\begin{dcases} \widetilde h(t), & t\in \overline{g({\mathbb R})},\\ h(c_n)+\frac{h(d_n)-h(c_n)}{d_n-c_n}(t-c_n), & t\in(c_n,d_n).
	\end{dcases}\] 
	$\widetilde f$ is clearly uniformly continuous by construction (as it is defined by saving the `gaps' in the graph of $\widetilde h$ with straight lines joining the endpoints) and $\widetilde f\circ g=f$.
\end{proof}
\begin{rem}The construction of $\widetilde f$ from $\widetilde h$ could also be carried out by adequate versions of the Tietze extension theorem for uniformly continuous functions on metric spaces \cite{Mandelkern1990}.
\end{rem}
\begin{rem}\label{remres}The domain of $f$ in the statement of Theorem~\ref{thmcomp} can be restricted to an interval $[a,b]$. Just extend $f$ constantly for values $t<a$ and $t>b$ and apply Theorem~\ref{thmcomp}. Then $\widetilde f\circ g|_{[a,b]}=f$.
\end{rem}
\begin{rem}
	We can present the Theorem~\ref{thmcomp} in terms of locally bounded locally uniformly $g$-continuous functions, that meaning that we can change the conditions in the statement to $f|_{[a,b]}\in\operatorname{BUC}_g([a,b],{\mathbb F})$ for every interval $[a,b]\subset {\mathbb R}$ and $\widetilde f\in\operatorname{C}({\mathbb R},{\mathbb F})$. We just have to apply Remark~\ref{remres} to ever increasing sequence of intervals $[a_n,b_n]$ to obtain extensions $\widetilde f_n$. Since the proof is constructive, the functions $\widetilde f_n|_{[a_n,b_n]}=\widetilde f_m$ for $[a_n,b_n]\subset[a_m,b_m]$ and we can thus define $\widetilde f$ on all of ${\mathbb R}$ satisfying the desired properties.

	Observe that, contrary to the classical case, not all $g$-continuous functions are locally bounded locally uniformly $g$-continuous functions \cite[Example 3.3]{FriLo17}. On the other hand, by Corollary~\ref{unif_reglada}, a function is locally bounded locally uniformly $g$-continuous function if and only if it is $g$-continuous and regulated.
\end{rem}

Let us denote by ${\mathbb F}[g]$ the ring of ${\mathbb F}$ polynomials on the variable $g$, that is, the set of all functions of the form $\sum_{k=0}^na_kg^k$ where $a_k\in {\mathbb F}$, $k=1,\dots, n$ and $n\in{\mathbb N}$. For the following result we will restrict the functions in ${\mathbb F}[g]$ to $[a,b]$.
\begin{thm}\label{WAP}
	${\mathbb F}[g]$ is dense in $\operatorname{BUC}_g([a,b],{\mathbb F})$.
\end{thm}
\begin{proof}
	Let $\varepsilon\in{\mathbb R}^+$, $f\in \operatorname{BUC}_g([a,b],{\mathbb F})$. By Theorem~\ref{thmcomp} and Remark~\ref{remres} $f=\widetilde f\circ g|_{[a,b]}$ for some uniformly continuous function $\widetilde f:{\mathbb R}\to{\mathbb R}$. Let $[c,d]\subset{\mathbb R}$ be such that $g([a,b])\subset[c,d]$. By the Weierstrass Approximation Theorem, there exists a polynomial $p\in{\mathbb F}[t]$ such that $\left\lVert \widetilde f|_{[c,d]}-p|_{[c,d]}\right\rVert_\infty<\varepsilon$. $p\circ g\in{\mathbb F}[g]$ and \[ \left\lVert f-p\circ g|_{[a,b]}\right\rVert_\infty =\left\lVert \widetilde f\circ g-p\circ g|_{[a,b]}\right\rVert_\infty\leqslant \left\lVert \widetilde f|_{[c,d]}-p|_{[c,d]}\right\rVert_\infty<\varepsilon.\qedhere\] 
\end{proof}

\begin{rem}The density of ${\mathbb R}[g]$ in $\operatorname{BUC}([a,b],{\mathbb F})$ comes with the important caveat that the functions in ${\mathbb R}[g]$ are not $g$-smooth due to the product rule \cite{Fernandez2022a}. Even so, as ${\mathbb Q}[g]$ is dense in ${\mathbb R}[g]$, we can conclude, as an straightforward consequence, that $\operatorname{BUC}([a,b],{\mathbb F})$ is a separable Banach space.\end{rem}

\section{Compactness in the space $\operatorname L_g^p(\mathbb{R},{\mathbb F})$}
In this section we prove an extension of Kolmogorov's compactness theorem for $L^p$ spaces. Before recalling the classical results --Theorems~\ref{k_r_2} and~\ref{k_r_1}), let us start with some technical results from \cite{Fernandez2022} that allow us to decompose Stieltjes integrals.
\begin{dfn}[Continuous and jump parts]\label{g_b_g_c}
	For a derivator $g$, we define $g^{B}:\mathbb{R}\to\mathbb{R}$ as
	\[ g^{B}(t):= \begin{dcases}\sum_{s \in[0, t) \cap D_{g}} \Delta g(s), & t>0, \\ -\sum_{s \in[t, 0) \cap D_{g}} \Delta g(s), & t \leqslant 0 ,\end{dcases}\] 
	which is left continuous and increasing. So is the mapping $g^{C}:\mathbb{R}\to\mathbb{R}$ given by
	\[ g^{C}(t):=g(t)-g^{B}(t).\] 
	We refer to $g^{C}$ and $g^{B}$ as the \emph{continuous} and \emph{jump} \emph{parts} of $g$, respectively. 
\end{dfn}
\begin{lem}[{\cite[Lemma 2.2]{Fernandez2022}}]
	$g$, $g^{C}$ and $g^{B}$ satisfy the following properties:
	\begin{enumerate}
		\item Given $A\in\mathcal{M}_g$, there exists $B\in G_{\delta}$ and $C\in{M}_g$ such that $A,C\subset B$, $\mu_g(C)=0$ and $A=B\backslash C$. 
		\item Given $A\in\mathcal{M}_g$, there exists $B\in F_{\sigma}$ and $C\in{M}_g$ such that $B\cap C=\emptyset$, $\mu_g(C)=0$ and $A=B\cup C$. 
		\item $\mathcal{M}_{g}\subset\mathcal{M}_{g^{C}}$.
		\item $\mathcal{M}_{g^{B}}=\mathcal{P}(\mathbb{R})$.
	\end{enumerate}
\end{lem}
Now we define a function similar to $\sigma$ occurring in the proof of Theorem~\ref{WAP}.

\begin{dfn}[Pseudoinverse of $g^C$]
	For a derivator $g:\mathbb{R}\to\mathbb{R}$ whose continuous part is not constant on intervals of the form $(-\infty,t]$, the \emph{pseudoinverse} of $g^{C}$ is the function 
	\[ \gamma: x \in g^{C}(\mathbb{R}) \rightarrow \gamma(x)=\min \left\{t \in\mathbb{R}: g^{C}(t)=x\right\}.\] 
\end{dfn}
\begin{pro}[{\cite[Proposition 2.6]{Fernandez2022}}]\label{pseudo_prop}
	The pseudoinverse of $g^{C}$ satisfies the following properties:
	\begin{itemize}
		\item For every $x\in g^{C}(\mathbb{R})$, $g^{C}(\gamma(x))=x$.
		\item For every $t\in\mathbb{R}$, $\gamma(g^{C}(t))\leqslant t$.
		\item For every $t\in\mathbb{R}\backslash(C_g\cup N_{g^{C}}^{+})$, $\gamma(g^{C}(t))=t$.
		\item $\gamma$ is strictly increasing.
		\item $\gamma$ is left continuous, and continuous in $x\in g^{C}(\mathbb{R})\backslash g^{C}(C_g)$.
	\end{itemize}
\end{pro}

\begin{lem}\label{descomposicion_1}
	For any $f\in{L}_{g}^{1}(\mathbb{R},{\mathbb F})$,
	\[ \int_{\mathbb{R}} f \mathrm{~d} \mu_{g}=\int_{\mathbb{R}} f \mathrm{~d} \mu_{g^{C}}+\sum_{s \in D_{g}} f(s) \Delta g(s).\] 
\end{lem}
\begin{proof}
	In \cite[Lemma 2.3]{Fernandez2022}, this result is proven for a bounded interval $[a,b)$. Applying this case and the dominated convergence theorem to $f_n=\chi_{[-n,n)}f$, the result for $\mathbb{R}$ follows.
\end{proof}
\begin{cor}
	For $E\in\mathcal{M}_g$,
	\[ \mu_{g}(E)=\mu_{g C}(E)+\sum_{s \in E \cap D_{g}} \Delta g(s).\] 
\end{cor}

The following result was also presented originally for a bounded interval.
\begin{pro}[{\cite[Proposition 2.7]{Fernandez2022}}]\label{morf_med}
	Given a derivator $g:\mathbb{R}\to\mathbb{R}$ and the Lebesgue $\sigma$-algebra $\mathcal{L}$:
	\begin{itemize}
		\item The continuous part $g^{C}:\left(\mathbb{R}, \mathcal{M}_{g}\right) \rightarrow\left(g^{C}(\mathbb{R}), \mathcal{L}\right)$ is a measurable morphism. 
		\item The pseudoinverse of the continuous part $\gamma:\left(g^{C}(\mathbb{R}), \mathcal{L}\right) \rightarrow\left( \mathbb{R}, \mathcal{M}_{g}\right)$ is a measurable morphism.
	\end{itemize}
\end{pro}
\begin{cor}\label{descomp_medida}
	Given $f\in{L}_g^{1}(\mathbb{R},{\mathbb F})$,
	\[ \int_{\mathbb{R}} f \mathrm{~d} \mu_{g}=\int_{g^{c}(\mathbb{R})} f \circ \gamma \mathrm{~d} \mu+\sum_{s \in D_{g}} f(s) \Delta g(s),\] 
	where $\mu$ is the Lebesgue measure.
\end{cor}

The proof of the following classical results regarding compactness in $L^{p}$ spaces can be found in \cite[Theorems 4,5]{HancheOlsen2010}.
\begin{thm}\label{k_r_2}
	If $1\leqslant p<\infty$, a subset of $\ell^{p}$ is totally bounded if and only if
	\begin{enumerate}
		\item it is pointwise bounded, 
		\item for every $\varepsilon>0$ there is an $n$ such that for any $x$ in the set,
		\[ \sum_{k>n}\left\lvert x_{k}\right\rvert^{p}<\varepsilon^{p}.\] 
	\end{enumerate}
\end{thm}
\begin{thm}[Kolmogorov-Riesz compactness theorem]\label{k_r_1}
	Suppose that $1\leqslant p<\infty$. A subset $\mathcal{F}$ of $L^{p}(\mathbb{R})$ is totally bounded if and only if
	\begin{enumerate}
		\item for every $\varepsilon>0$ there exists $R$ such that, for every $f\in\mathcal{F}$,
		\[ \int_{\left\lvert x\right\rvert>R}\left\lvert f(x)\right\rvert^{p}\operatorname{d} x<\varepsilon^{p},\] 
		\item for every $\varepsilon>0$ there exists $\rho>0$ such that for every $f\in\mathcal{F}$ and $h\in\mathbb{R}$ with $\left\lvert h\right\rvert<\rho$,
		\[ \int_{\mathbb{R}}\left\lvert f(x+h)-f(x)\right\rvert^{p}\operatorname{d} x<\varepsilon^{p}.\]  
	\end{enumerate} 
\end{thm}
\begin{rem}\label{redundante}
	The Theorem~\ref{k_r_1} is often stated with the additional condition that $\mathcal{F}$ is a bounded subset of $L^{p}(\mathbb{R},{\mathbb F})$, but it is possible to prove that this is implied by conditions $1$ and $2$ --see \cite{HancheOlsen2019}. 
\end{rem}

Before moving on to the main result of the section, we need to make one last consideration. The pseudoinverse of $g^C$, $\gamma$, is defined on the set $g^{C}(\mathbb{R})$. Consequently, for a family $\mathcal{G}\subset L_{g}^{p}(\mathbb{R},{\mathbb F})$, $\mathcal{G}\circ\gamma:=\{f\circ\gamma:f\in\mathcal{G}\}$ is a subset of $L^{p}(g^{C}(\mathbb{R}),{\mathbb F})$. Since we can embed this space in $L^{p}(\mathbb{R})$ (defining its functions as zero outside of $g^{C}(\mathbb{R})$), we will treat $\mathcal{G}\circ\gamma$ as a subset of $L^{p}(\mathbb{R},{\mathbb F})$. 	
\begin{thm}\label{k-r-g}
	Suppose that $1\leqslant p<\infty$ and $D_g=\{d_1,d_2,...\}$. A subset $\mathcal{G}$ of $L_{g}^{p}(\mathbb{R},{\mathbb F})$ is totally bounded if and only if
	\begin{enumerate}
		\item $\{f(d_n):f\in\mathcal{G}\}\subset\mathbb{R}$ is bounded for each $d_n\in D_g$,
		\item for every $\varepsilon>0$ there exists $n$ such that for any $f\in\mathcal{G}$,
		\[ \sum_{k>n}\left\lvert f(d_k)\right\rvert^{p}\Delta g(d_k)<\varepsilon^p,\] 
		\item for every $\varepsilon>0$ there is some $R$ such that, for any $f\in\mathcal{G}$,
		\[ \int_{\left\lvert x\right\rvert>R}\left\lvert f \circ \gamma(x)\right\rvert^{p} \operatorname{d} x<\varepsilon^{p},\] 
		\item for every $\varepsilon>0$ there exists $\rho>0$ such that for every $f\in\mathcal{G}$ and $h\in\mathbb{R}$ with $\left\lvert h\right\rvert<\rho$,
		\[ \int_{\mathbb{R}}\left\lvert f\circ\gamma(x+h)-f\circ\gamma(x)\right\rvert^{p} \operatorname{d} x<\varepsilon^{p}.\] 
	\end{enumerate}
\end{thm}
\begin{proof}
	If $f,h\in L_g^{p}(\mathbb{R},{\mathbb F})$, thanks to Corollary~\ref{descomp_medida}, we have that
	\begin{equation}\label{eq_k_r}
		\begin{aligned}
			\left\lVert f-h\right\rVert_{L_g^{p}}^{p}=\int_{\mathbb{R}} \left\lvert f-h\right\rvert^{p} \mathrm{~d} \mu_{g} =&\int_{g^{c}(\mathbb{R})} \left\lvert f \circ \gamma-h \circ \gamma\right\rvert^{p} \mathrm{~d} \mu+\sum_{s \in D_{g}} \left\lvert f(s)-h(s)\right\rvert^{p} \Delta g(s)\\ =&\left\lVert f \circ \gamma-h \circ \gamma\right\rVert_{L^{p}}^{p}+\left\lVert ((f(d_n)-h(d_n))\Delta g(d_n)^{1/p})_{n\in\mathbb{N}}\right\rVert_{\ell^{p}}^{p}.
		\end{aligned}
	\end{equation}\par
	Suppose that $\mathcal{G}\subset L_g^{p}(\mathbb{R},{\mathbb F})$ is totally bounded. With the previous equality in mind, it is clear that an $\varepsilon$-net of $\mathcal{G}$ induces an $\varepsilon$-net on the sets $\{f\circ\gamma:f\in\mathcal{G}\}\subset L^{p}(\mathbb{R},{\mathbb F})$ and $\{(f(d_n)\Delta g(d_n)^{1/p})_{n\in\mathbb{N}}:f\in\mathcal{G}\}\subset \ell^{p}$, which are then totally bounded. Using Theorems~\ref{k_r_2} and~\ref{k_r_1}, we conclude that $\mathcal{G}$ satisfies conditions 1-4 (for condition 1, note that for each fixed $d_n$, $\{f(d_n)\Delta g(d_n)^{1/p}:f\in\mathcal{G}\}$ is bounded if and only if $\{f(d_n):f\in\mathcal{G}\}$ is).\par
	Let us finally check that if all of the conditions are satisfied, every sequence in $\mathcal{G}$ has a Cauchy subsequence, which is equivalent to $\mathcal{G}$ being totally bounded. Take $\{f_n\}_{n\in\mathbb{N}}\subset\mathcal{G}$. The sets $\{f\circ\gamma:f\in\mathcal{G}\}\subset L^{p}(\mathbb{R},{\mathbb F})$ and $\{(f(d_n)\Delta g(d_n)^{1/p})_{n\in\mathbb{N}}\}\subset \ell^{p}$ satisfy the hypotheses of Theorems~\ref{k_r_2} and~\ref{k_r_1}, respectively, so $\{f_n\circ\gamma\}_{n\in\mathbb{N}}$ admits a Cauchy subsequence $\{f_{n_k}\circ\gamma\}_{k\in\mathbb{N}}$ in $L^{p}(\mathbb{R},{\mathbb F})$ and $\{(f_{n_k}(d_j)\Delta g(d_j)^{1/p})_{j\in\mathbb{N}}\}_{k\in\mathbb{N}}$ also admits a Cauchy subsequence $\{(f_{n_{k_l}}(d_j)\Delta g(d_j)^{1/p})_{j\in\mathbb{N}}\}_{l\in\mathbb{N}}$ in $\ell^{p}$. Again thanks to~\eqref{eq_k_r}, $\{f_{n_{k_l}}\}_{l\in\mathbb{N}}$ is a Cauchy subsequence of $\{f_n\}_{n\in\mathbb{N}}$ in $L_g^{p}(\mathbb{R},{\mathbb F})$. 
\end{proof}
\begin{rem}
	As in Theorem~\ref{k_r_1}, conditions $3$ and $4$ imply that $\mathcal{G}\circ\gamma$ is a bounded subset of $L^p(\mathbb{R},{\mathbb F})$. 
\end{rem}

\section{An application to Stieltjes-Sobolev spaces}

In this section we define the Stieltjes-Sobolev spaces and generalize the classical results concerning continuous and compact inclusions of Sobolev space theory. In particular, these compactness results will be useful when it comes to studying nonlinear Stieltjes differential equations.


When integrating funtions, we will use the following notation for convenience:

\[ \int_x^y\varphi(s)\, \operatorname{d}\mu_g(s)=\left\{
\begin{aligned}
	&+ \int_{[x,y)} \varphi(s)\, \operatorname{d}\mu_g(s),\; y\geq x,\\
	& -\int_{[y,x)} \varphi(s)\, \operatorname{d}\mu_g(s),\; y\leqslant x.
\end{aligned}
\right.\] 

We start by defining the concept of Stieltjes-Sobolev space $W^{1,p}_g(I,{\mathbb F})$.

\begin{dfn}[The Stieltjes-Sobolev space $W^{1,p}_g(I,{\mathbb F})$]\label{dfnsss} Let $I\subset 
		\mathbb{R}$ be an interval which is closed from the left (in the case it is bounded from below) and open from the right, and $p\in[1,\infty]$. We define
		\begin{displaymath}
			W_g^{1,p}(I,{\mathbb F})=\left\{ u\in L^p_g(I):\: \exists \widetilde{u}\in L_g^p(I)\text{ such that }
			u(y)-u(x)=\int_x^y \widetilde{u}(s)\, \operatorname{d}\mu_g(s),\; \forall x,y\in \overline{I}
			\right\}.
		\end{displaymath}

		We endow $W_g^{1,p}(I,{\mathbb F})$ with the norm
		\begin{displaymath}
			\left\lVert u\right\rVert_{W_g^{1,p}(I)}:=\left\lVert u\right\rVert_{L^p_g(I)}+\left\lVert \widetilde{u}\right\rVert_{L^p_g(I)},
\end{displaymath}
and, we will write, in general, $W_g^{1,p}(I)\equiv W_g^{1,p}(I,{\mathbb F})$ and follow the same convention for the rest of spaces.
\end{dfn}

\begin{rem} Observe that
		 $L^p_g(I)\subset L_{g,\rm{loc}}^1(I)$, so Definition~\ref{dfnsss} makes sense, and 
		 $W_g^{1,p}(I)\subset \operatorname{C}_g(\overline{I})$, for any $p\in [1,\infty]$. 
		\end{rem}

		\begin{rem} 
 If $I=[a,b)$ (where $b$ can take the value $\infty$), we can write
		\begin{displaymath}
				W_g^{1,p}([a,b))=\left\{ u\in L^p_g([a,b)):\: \exists \widetilde{u}\in L_g^p([a,b)),\ 
				u(x)=u(a)+\int_{[a,x)} \widetilde{u}(s)\, \operatorname{d}\mu_g(s),\; \forall x \in [a,b]
				\right\}.
		\end{displaymath}
		In the case $b=\infty$, we take $x\in[a,b)$.
\end{rem}

\begin{rem} In the particular case $a,b\in \mathbb{R}$, $a<b$, since $\mu_g([a,b))=g(b)-g(a)<\infty$, 
it is clear that ${L}_g^p([a,b))\subset 
{L}_g^1([a,b))$, for every $ p\in[1,\infty]$. Thus, thanks to \cite[Theorem 2.4]{LoRo14}, 
there exists a $g$-measurable set $N\subset [a,b]$, with $\mu_g(N)=0$, such that 
$u_g'(t)=\widetilde{u}(t)$, for all $t\in [a,b]\setminus N$. Moreover, thanks to 
\cite[Proposition 5.2]{LoRo14}, $W_g^{1,p}([a,b))\subset \operatorname{AC}_g([a,b])$, 
$p\in[1,\infty]$. Therefore, the elements of the space $W^{1,p}_g([a,b))$ are 
$g$-absolutely continuous functions whose $g$-derivatives are in the space 
${L}_g^p([a,b))$.

On the other hand --see \cite[Proposition5.5]{FriLo17}, 
$W_g^{1,p}([a,b))\subset \operatorname{AC}_g([a,b])\subset \operatorname{BC}_g([a,b])$. In the 
next lemma we will further show that the embedding $W_g^{1,p}([a,b))\subset 
\operatorname{BC}_g([a,b])$ is continuous.
\end{rem}

\begin{lem} \label{embc1} Let us consider $a,b\in \mathbb{R}$ such that $a<b$. The embedding $W_g^{1,p}([a,b))\subset 
\operatorname{BC}_g([a,b])$ is continuous for every $ p\in[1,\infty]$.
\end{lem}

\begin{proof} Given $u\in W_g^{1,p}([a,b))$, we define:
\begin{displaymath}
v(t)=\int_0^t \widetilde{u}(s)\, \operatorname{d}\mu_g(s).
\end{displaymath}
We have that $v\in \operatorname{AC}_g([a,b])$. Moreover, if $p=1$, 
\begin{displaymath}
\left\lvert v(t)\right\rvert\leqslant \int_{[a,t)}\left\lvert \widetilde{u}(s)\right\rvert\, \operatorname{d}\mu_g(s)
\leqslant \left\lVert \widetilde{u}\right\rVert_{{L}_g^1([a,b))},\;
\forall t\in [a,b].
\end{displaymath}
If $1<p<\infty$,
\begin{displaymath}
\left\lvert v(t)\right\rvert\leqslant \int_{[a,b)} \left\lvert \widetilde{u}(s)\right\rvert\, \operatorname{d}\mu_g(s)
\leqslant \mu_g([a,b))^{p/(p-1)}\, \left\lVert \widetilde{u}\right\rVert_{{L}_g^p([a,b))},\;
\forall t\in [a,b).
\end{displaymath}
If $p=\infty$,
\begin{displaymath}
\left\lvert v(t)\right\rvert\leqslant \int_{[a,b)} \left\lvert \widetilde{u}(s)\right\rvert\leqslant \mu_g([a,b)) \, 
\left\lVert \widetilde{u}\right\rVert_{{L}_g^{\infty}([a,b))},\;
\forall t\in [a,b).
\end{displaymath}
Therefore, given $1\leqslant p\leqslant \infty$, there exists a positive 
constant $C\equiv C(p,\mu_g([a,b)))>0$ such that,
\begin{displaymath}
\left\lvert v(t)\right\rvert\leqslant C\, \left\lVert \widetilde{u}\right\rVert_{{L}_g^p([a,b))}, 
\; \forall t \in [a,b].
\end{displaymath}
Thus, there is a constant, which we will continue to denote in the same way, such that
\begin{displaymath}
\left\lVert v\right\rVert_{{L}_g^p([a,b))} \leqslant 
C\, \left\lVert \widetilde{u}\right\rVert_{{L}_g^p([a,b))}.
\end{displaymath}

Now, we have that $u(t)=v(t)+u(a)$, in particular $u(a)=u(t)-v(t)$, so there exists another 
positive constant, that we will continue to denote in the same way, such that
\begin{displaymath}
\left\lvert u(a)\right\rvert\leqslant C\,\left(\left\lVert u\right\rVert_{{L}_g^p([a,b))}+
\left\lVert \widetilde{u}\right\rVert_{{L}_g^p([a,b))} \right).
\end{displaymath}
Finally,
\begin{align*}
\left\lVert u\right\rVert_{0}=&\sup_{t\in [a,b]} \left\lvert u(t)\right\rvert 
=\sup_{t\in [a,b]} \left\lvert  v(t)+u(a)\right\rvert
\leqslant \sup_{t\in [a,b]} \left\lvert  v(t)\right\rvert+\left\lvert u(a)\right\rvert\\
\leqslant & C\, \left\lVert \widetilde{u}\right\rVert_{{L}_g^p([a,b))}+
C\,\left(\left\lVert u\right\rVert_{{L}_g^p([a,b))}+
\left\lVert \widetilde{u}\right\rVert_{{L}_g^p([a,b))} \right)\\
\leqslant & 2\, C\, \left(\left\lVert u\right\rVert_{{L}_g^p([a,b))}+
\left\lVert \widetilde{u}\right\rVert_{{L}_g^p([a,b))} \right).
\qedhere
\end{align*}
\end{proof}

\begin{thm} Let $I\subset \mathbb{R}$ be an interval which is closed from the left (in the case it is bounded from below) and open from the right. We have that $W_g^{1,p}(I)$ is a Banach space 
for all $p\in [1,\infty]$. Furthermore, if $p\in (1,\infty)$, $W_g^{1,p}(I)$ is reflexive and, for $p\in [1,\infty)$, $W_g^{1,p}(I)$ is 
separable.
\end{thm}
\begin{proof} Let us consider a Cauchy sequence $\{u_n\}_{n \in \mathbb{N}}\subset W_g^{1,p}(I)$. We have, thanks to the definition of the norm in $W_g^{1,p}(I)$ and Lemma~\ref{embc1}, that
$\{u_n\}_{n \in \mathbb{N}}$ is a Cauchy sequence in ${L}_g^p(I)$ 
and $\{\widetilde{u}_n\}_{n \in \mathbb{N}}$ is a Cauchy sequence in ${L}_g^p(I)$. 
Moreover, due to Lemma~\ref{embc1}, $\{u_n\}_{n \in \mathbb{N}}$ is a Cauchy sequence 
in $\operatorname{BC}_g([a,b])$, for all $a,\,b\in \mathbb{R}$, $a<b$.

Thanks to the completeness of ${L}_g^p(I)$, there exist elements 
$u,\,\widetilde{u}\in {L}_g^p(I)$ such that
\begin{displaymath}
\begin{aligned}
\lim_{n \to \infty} \left\lVert u_n-u\right\rVert_{{L}_g^p(I)} =&0, \\
\lim_{n \to \infty} \left\lVert \widetilde{u}_n-\widetilde{u}\right\rVert_{{L}_g^p(I)} =&0.
\end{aligned}
\end{displaymath}	

Now, given $x,\, y \in \overline{I}$, $x<y$, 
\begin{equation}\label{eqconv}
\left\lvert  \int_{[x,y)} \left(\widetilde{u}_n(s)-\widetilde{u}(s)\right)\, \operatorname{d}\mu_g(s)\right\rvert
\leqslant \int_{[x,y)} \left\lvert \widetilde{u}_n(s)-\widetilde{u}(s)\right\rvert\, \operatorname{d}\mu_g(s)
\leqslant C\, \left\lVert \widetilde{u}_n-\widetilde{u}\right\rVert_{{L}_g^p(I)},
\end{equation}
where $C\equiv C(p,\mu_g([x,y)))>0$ is a positive constant. Thus, because of inequality~\eqref{eqconv} and the fact that $\lim_{n \to \infty} \left\lVert \widetilde{u}_n-\widetilde{u}\right\rVert_{{L}_g^p(I)} =0$, we have that, for $x,\, y \in \overline{I}$, $x<y$, 
\begin{equation} \label{eq:cl1}
\lim_{n \to \infty} \int_{[x,y)} \widetilde{u}_n(s) \operatorname{d}\mu_g(s) = \int_{[x,y)} \widetilde{u}(s)\, \operatorname{d}\mu_g(s), \; 
\forall x,\,y\in \mathbb{R}.
\end{equation}

On the other hand, since $\operatorname{BC}_g([a,b])$ is a Banach space \cite[Theorem 3.4]{FriLo17}
for all $a,\,b\in \mathbb{R}$, $a<b$, and the embedding $W_g^{1,p}([a,b))\subset 
\operatorname{BC}_g([a,b])$ is continuous, for all $p\in[1,\infty]$, we have that $u|_{[a,b]}\in \mathcal{C}_g([a,b])$ and 
\begin{displaymath}
\lim_{n\to \infty} u_n(t)=u(t),\; \forall t\in [a,b].
\end{displaymath}
Since $[a,b]$ was fixed arbitrarily, $u\in \mathcal{C}_g(I)$ and 
\begin{equation} \label{eq:cl2}
	\lim_{n\to \infty} u_n(t)=u(t),\; \forall t\in \overline I.
\end{equation}
Finally, thanks to~\eqref{eq:cl1} and~\eqref{eq:cl2}, we can pass to the limit in the following expression
\begin{displaymath}
u_n(y)=u_n(x)+\int_x^y \widetilde{u}_n(s)\, \operatorname{d}\mu_g(s),\; x,\,y\in \overline{I},\; n\in \mathbb{N},
\end{displaymath}
and we obtain that
\begin{displaymath}
u(y)=u(x)+\int_x^y \widetilde{u}(s)\, \operatorname{d}\mu_g(s),\; 
\forall x,\, y \in \overline{I}.
\end{displaymath}
Thus, $u_n$ converges to $u$ in $W_g^{1,p}(I)$.

Now, given $p\in (1,\infty)$, we have by \cite[Theorem B.92]{Leoni2017} that 
$L^p_g(I)$ is reflexive. Thus the product space $\widetilde{W}=L^p_g(I)\times L^p_g(I)$ is 
also reflexive. The operator $T:u\in W_g^{1,p}(I)\rightarrow T(u)=(u,\widetilde{u})\in \widetilde{W}$ is an isometry 
from $W_g^{1,p}(I)$ to $\widetilde{W}$. Since $W_g^{1,p}(I)$ is a Banach space, $T(W_g^{1,p}(I))$ is a 
closed subspace of $\widetilde{W}$. It now follows that $T(W_g^{1,p}(I))$ is reflexive, so $W_g^{1,p}(I)$ is also reflexive.

Finally, since $L^p_g(I)$ is separable \cite[Theorem 4.13]{Brezis2011} for $p\in [1,\infty)$, we have that $\widetilde{W}=L^p_g(I)\times L^p_g(I)$ 
is also separable and, therefore, $T(W_g^{1,p}(I))$ is separable, implying that $W_g^{1,p}(I)$ is separable.
\end{proof}

\subsection{Compact embedding into $\operatorname{BC}_g([a,b])$}

One of the classical results in Sobolev spaces is the compact embedding of $W^{1,p}(a,b)$ into $\operatorname{C}([a,b])$ 
for $p\in(1,\infty]$ \cite[Theorem 8.8]{Brezis2011}. In the following Theorem we generalize the previous result to the 
case of Stieltjes-Sobolev spaces, that is, we will prove that the embedding $W_g^{1,p}([a,b))\subset \operatorname{BC}_g([a,b])$ 
is compact, for all $p\in(1,\infty]$.

\begin{thm}[Compactness of $W_g^{1,p}([a,b))$ into $\operatorname{BC}_g{([a,b])}$] \label{cont3} The embedding of $W_g^{1,p}([a,b))$ into $\operatorname{BC}_g{([a,b])}$ is compact for all $p\in(1,\infty]$.
\end{thm}
\begin{proof}Let us consider $a,b\in \mathbb{R}$ such that $a<b$, $p\in(1,\infty]$ and define
	\begin{displaymath}
		\mathcal{H}:=\{u \in W_g^{1,p}([a,b)):\; \left\lVert u\right\rVert_{W^{1,p}_g(a,b)}\leqslant 1\}.
	\end{displaymath}
To show that the embedding of $W_g^{1,p}([a,b))$ into $\operatorname{BC}_g{([a,b])}$ is compact it is enough to show that $\mathcal{H}$ is a relatively compact subset of $\operatorname{BC}_g([a,b])$. Now, the proof is a direct consequence of Corollary~\ref{cor:precompact}. Indeed, 
 first, the set $\{u(t):\; u \in \mathcal{H}\}$ is bounded for all $t\in [a,b]$. Given an 
element $u\in \mathcal{H}$,
\begin{displaymath}
\left\lvert u(t)\right\rvert\leqslant \left\lVert u\right\rVert_0 \leqslant \left\lVert u\right\rVert_{W^{1,p}_g([a,b])}\leqslant 1.
\end{displaymath}

Second, $\mathcal{H}$ is uniformly $g$-equicontinuous. Given $t,\, s\in [a,b]$ with $t\leqslant s$ and $u \in \mathcal{H}$,
\begin{displaymath}
\left\lvert u(s)-u(t)\right\rvert\leqslant \left\lvert \int_{[t,s)} u_g'(u)\, \operatorname{d}\mu_g(u)\right\rvert \leqslant \left\lVert u_g'\right\rVert_{{L}_g^p([a,b))} \left[ \int_{[t,s)} 1\, \operatorname{d}\mu_g(s)\right|^{1/p'}
\leqslant \left[g(s)-g(t)\right]^{1/p'},
\end{displaymath}
where $p'=1$ if $p=\infty$ and $1/p+1/p'=1$ otherwise. Therefore, given $\varepsilon>0$, we can take $\delta=\varepsilon^{p'}$ and, then, 
\begin{displaymath}
\left\lvert u(s)-u(t)\right\rvert\leqslant \left[g(s)-g(t)\right]^{1/p'} < \delta^{1/p'}=\varepsilon, 
\end{displaymath}
for every $u\in \mathcal{H}$ and every $t,s\in [a,b]$ such that $\left\lvert g(t)-g(s)\right\rvert<\delta$.
\end{proof}

\begin{rem} Observe that, thanks to Proposition~\ref{equi_unif_equi}, we can replace the condition in the definition of $\mathcal{H}$ in the proof of Theorem~\ref{cont3} by asking for the family ${\mathcal H}$ to have
 uniform right-hand side limits in $[a,b)\cap D_g$. To prove this equivalent 
condition it suffices to observe that given an element $t_0\in [a,b)\cap D_g$, we have that
\begin{displaymath}
\left\lvert u(t)-u(t_0^+)\right\rvert=\left\lvert \int_{(t_0,t)} u_g'(s)\, \operatorname{d}\mu_g(s) \right\rvert \leqslant 
\left[ g(t)-g(t_0^+)\right]^{1/p'},\; \forall t\in (t_0,b],\; \forall u \in \mathcal{H}.
\end{displaymath}
Now, by definition of the right-hand side limit at $t_0$, there exists an element $\delta>0$, such that 
if $t\in (t_0,t_0+\delta)$, then $g( t) -g(t_0^ +)<\varepsilon^{p'}$. Then, taking $\delta=\varepsilon^{p'}$ 
we have the desired result.
\end{rem}

\subsection{Compact embedding into $L_g^q([a,b))$}

In this section we generalize the compact embedding of $W^{1,1}(a,b)$ into 
$L^q(a,b)$ for all $q\in [1,\infty)$ \cite[Theorem 8.8]{Brezis2011} to the general case of 
Stieltjes-Sobolev spaces by using Theorem~\ref{k-r-g}. We must emphasize that the classical proof of this result 
goes through the construction of an extension operator $P:W^{1,1}(a,b)\rightarrow W^{1,1}(\mathbb{R})$, for 
which reflexion and prolongation techniques are used. In the case of the Stieltjes derivative, the reflexion techniques 
are not directly applicable, so we propose an alternative construction of an extension operator based on concatenation 
with $g$-exponential functions. To do this, it will be necessary to extend the definition the 
$g$-exponential function in \cite[Definition 4.5]{Fernandez2022} to the whole real line.

\begin{dfn}[global $g$-exponential]\label{gge}
	Let $\alpha\in 
	\mathbb{R}$ and 
	$\lambda \in L^1_{g,\rm{loc}}(\mathbb{R})$ such that $1+\lambda(t)\, \Delta g(t) \neq 0$, 
	$\forall t\in D_g$. Define, for $t\geqslant \alpha$, 	\begin{displaymath}
		\exp_g(\lambda;\alpha,t)=\exp\left(\int_{[\alpha,t)}\widetilde{\lambda}(s)\, \operatorname{d}\mu_g(s)\right),\; \forall \alpha, t\in \mathbb{R},
	\end{displaymath}
	where 
	\begin{displaymath}
		\widetilde{\lambda}:t\in \mathbb{R} \rightarrow \widetilde{\lambda}(t)=
		\begin{dcases}
			 \lambda(t), & t\in \mathbb{R}\setminus D_g, \\
		\frac{\log\left(1+\lambda(t)\, \Delta g(t)\right)}{\Delta g(t)}, & t\in D_g.
		\end{dcases}
	\end{displaymath}
and define, for $t<\alpha$,
		\begin{displaymath}
		\exp_g(\lambda;\alpha,t):=
			\exp_g(\lambda;t,\alpha)^{-1}=\exp_g(q(\lambda); t, \alpha),
	\end{displaymath}
	where
	\begin{displaymath}
		q(\lambda)(t)=-\frac{\lambda(t)}{1+\lambda(t)\,\Delta g(t)}.
	\end{displaymath}
\end{dfn}

\begin{rem} Observe that, thanks to \cite[Proposition~4.6]{Fernandez2022}, given $\lambda \in L^1_{g,\rm{loc}}(\mathbb{R})$ such 
that $1+\lambda(t)\, \Delta g(t) \neq 0$, $\forall t\in D_g$, $q(\lambda)\in L^1_{g,\rm{loc}}(\mathbb{R})$ and 
$1+q(\lambda)(t)\, \Delta g(t) \neq 0$, $\forall t\in D_g$. On the other hand, If $\lambda \in L^1_{g,\rm{loc}}(\mathbb{R})$ 
is such that $1+\lambda(t)\Delta g(t)\neq 0$, $\forall t\in \mathbb{R}\cap D_g$, we have by \cite[Remark 4.4]{Fernandez2022} that 
$\widetilde{\lambda}\in L^1_{g,\rm{loc}}(\mathbb{R})$. Therefore, the global $g$-exponential is well defined.
\end{rem}

For the next result, we denote the set of locally bounded functions defined on ${\mathbb R}$ as $\operatorname{B}_{\rm{loc}}(\mathbb{R})$.
\begin{lem}
	If $\lambda \in \operatorname{C}_{g}(\mathbb{R})\cap 
	\operatorname{B}_{\rm{loc}}(\mathbb{R})$, then
	\begin{displaymath}
		\left(\exp_g(\lambda;\alpha, \cdot )\right)'_g(t)=\lambda(t)\, \exp_g(\lambda;\alpha,t), \;
		\forall t\geqslant \alpha.
	\end{displaymath}
\end{lem}


\begin{lem} \label{lemgexpg} Let $\alpha\in 
	\mathbb{R}$ and 
	$\lambda \in L^1_{g,\rm{loc}}(\mathbb{R})$ such that $1+\lambda(t)\, \Delta g(t) \neq 0$, 
	$\forall t\in D_g$. We have that $\exp_g(\lambda;\alpha,\cdot)\in \operatorname{AC}_{g,\rm{loc}}(\mathbb{R})\subset \operatorname{C}_{g}(\mathbb{R})\cap 
	\operatorname{B}_{\rm{loc}}(\mathbb{R})$. Furthermore, given $x,y\in \mathbb{R}$,
	\begin{displaymath}
		\exp_g(\lambda;\alpha,y)-\exp_g(\lambda;\alpha,x)=\int_x^y \lambda(s)\, 
		\exp_g(\lambda;\alpha,s)\, \operatorname{d}\mu_g(s).
	\end{displaymath}
\end{lem}

\begin{proof} We will assume, without loss of generality, that $x\leqslant y$. We study different cases.

	If $y>x\geq \alpha$, the results follows from the definition of the $g$-exponential. On the other hand, for $x<y<\alpha$, given $\beta<x$, 
		\begin{align*}
				&\exp_g(\lambda;\alpha,y)-\exp_g(\lambda;\alpha,x)\\
				=& \exp\left(\int_y^{\alpha} 
				\widetilde{q(\lambda)}(s)\, \operatorname{d}\mu_g(s) \right)-
				\exp\left(\int_x^{\alpha} 
				\widetilde{q(\lambda)}(s)\, \operatorname{d}\mu_g(s) \right)
				\\
				=&
				\exp\left(\int_{\beta}^{\alpha}
				\widetilde{q(\lambda)}(s)\, \operatorname{d}\mu_g(s)-
				\int_{\beta}^y 
				\widetilde{q(\lambda)}(s)\, \operatorname{d}\mu_g(s) \right)\\
				&
				-\exp\left(\int_{\beta}^{\alpha}
				\widetilde{q(\lambda)}(s)\, \operatorname{d}\mu_g(s)-
				\int_{\beta}^x 
				\widetilde{q(\lambda)}(s)\, \operatorname{d}\mu_g(s) \right)
				\\
				=&
				\exp\left(\int_{\beta}^{\alpha}
				\widetilde{q(\lambda)}(s)\, \operatorname{d}\mu_g(s)\right)
				\left(\left[\exp\left(\int_{\beta}^y 
				\widetilde{q(\lambda)}(s)\, \operatorname{d}\mu_g(s) 
				\right) \right]^{-1} - 
				\left[\exp\left(\int_{\beta}^x 
				\widetilde{q(\lambda)}(s)\, \operatorname{d}\mu_g(s) 
				\right) \right]^{-1}
				\right)
				\\
			=	&
				\exp\left(\int_\beta^\alpha \widetilde{q(\lambda)}(s)\, \operatorname{d}\mu_g(s) \right) 
				\left(\left[\exp_g(q(\lambda);\beta,y)\right]^{-1}-\left[ \exp_g(q(\lambda);\beta,x) \right]^{-1} \right)\\
				=&\exp\left(\int_\beta^\alpha \widetilde{q(\lambda)}(s)\, \operatorname{d}\mu_g(s) \right) 
				\left(\exp_g(\lambda;\beta,y)-\exp_g(\lambda;\beta,x) \right)\\
				=&\exp\left(\int_\beta^\alpha \widetilde{q(\lambda)}(s)\, \operatorname{d}\mu_g(s) \right) 
				\int_x^y\lambda(s)\,\exp_g(\lambda;\beta,s)\, \operatorname{d}\mu_g(s)\\
				=& \exp\left(\int_\beta^\alpha \widetilde{q(\lambda)}(s)\, \operatorname{d}\mu_g(s) \right) 
				\int_x^y\lambda(s) \frac{1}{\exp_g(q(\lambda);\beta,s)}\, \operatorname{d}\mu_g(s)\\
				=&\int_x^y \lambda(s)\, \exp\left( \int_\beta^\alpha \widetilde{q(\lambda)}(u)\, \operatorname{d}\mu_g(u) -
				\int_\beta^s\widetilde{q(\lambda)}(u)\, \operatorname{d}\mu_g(u) \right)\, \operatorname{d}\mu_g(s)\\
				=&\int_x^y \lambda(s)\, \exp\left(\int_s^\alpha \widetilde{q(\lambda)}(u)\, \operatorname{d}\mu_g(u)\right)\, 
				\operatorname{d}\mu_g(s)=\int_x^y \lambda(s)\,\exp_g(\lambda;\alpha,s)\, \operatorname{d}\mu_g(s).\qedhere
		\end{align*}
	Finally, in the case $x\leqslant \alpha\leqslant y$, we can proceed by combining the two previous cases.\end{proof}

Let us study when the $g$-exponential function belongs to
$W_g^{1,p}([\alpha,\infty))$ and $W_g^{1,p}((-\infty,\alpha))$.

\begin{lem} \label{gexpgext1} Let $\lambda>0$ be such that $1-\lambda \Delta g(t)\in (0,1)$, $\forall t\in 
	[\alpha,\infty)\cap D_g$. We have that $\exp_g(-\lambda;\alpha,\cdot) 
	\in W_g^{1,p}([\alpha,\infty))$, for every $p\in [1,\infty]$. Furthermore,
	$\left\lVert u\right\rVert_{W_g^{1,p}([\alpha,\infty))}\leqslant C(\lambda)$, where $C(\lambda)>0$ is a constant that does not depend on $p$.
\end{lem}

\begin{proof} On the one hand, by the definition of the $g$-exponential 
 and the hypotheses on $\lambda$, we have that
	\begin{displaymath}
		\exp_g(-\lambda;\alpha,t)=\exp\left(
		-\lambda \mu_g([\alpha,t)\setminus D_g)+\sum_{s\in [\alpha,t)} \log\left(
		1-\lambda\, \Delta g(s)
		\right)
		\right)\in (0,1],\; \forall t\in [\alpha,\infty).
	\end{displaymath}
	On the other, by the properties of the $g$-exponential and for $\alpha \leqslant x<y$ --see Lemma~\ref{lemgexpg}, 
	we deduce that
	\begin{displaymath}
		\exp_g(-\lambda;\alpha,y)-\exp_g(-\lambda;\alpha,x)=
		-\lambda\int_x^y \exp_g(-\lambda;\alpha,s)\, \operatorname{d}\mu_g(s).
	\end{displaymath}
With these observations we can now prove that 
	$\exp_g(-\lambda;\alpha,\cdot) \in L_g^p([\alpha,\infty))$, 
	for every $p\in [1,\infty)$. We have two cases: either $p=\infty$, case for which the result is obvious since $\exp_g(-\lambda;\alpha,t)\in (0,1]$, 
	$\forall t\in [\alpha,\infty)$, or $p\in[1,\infty)$. In this case, the monotone convergence theorem,
\begin{displaymath}
\int_{[\alpha,\infty)} \left\lvert \exp_g(-\lambda;\alpha,s)\right\rvert^p\, \operatorname{d}\mu_g(s)=\lim_{M \to \infty} \int_{\alpha}^M \left\lvert \exp_g(-\lambda;\alpha,s)\right\rvert^p\, 
\operatorname{d}\mu_g(s)
\end{displaymath}
Now, given that $\exp_g(-\lambda;\alpha,t)\in (0,1]$ we have that
\begin{displaymath}
\exp_g(-\lambda;\alpha,t)^p\leqslant \exp_g(-\lambda;\alpha,t),\; \forall t\in [\alpha,\infty),\; 
\forall p\in [1,\infty),
\end{displaymath}
and we deduce that
		\begin{displaymath}
			\begin{aligned}
				\int_{\alpha}^{M} \exp_g(-\lambda;\alpha,s)^p\, \operatorname{d}\mu_g(s)\leqslant & 
				\int_{\alpha}^{M} \exp_g(-\lambda;\alpha,s)\, \operatorname{d}\mu_g(s) \\
				=&\frac{1}{-\lambda} \int_{\alpha}^{M} -\lambda\, \exp_g(-\lambda;\alpha,s)\, \operatorname{d}\mu_g(s)
				= \frac{1}{-\lambda} \left( \exp_g(-\lambda;\alpha,M)-1 \right)\\
				=&\frac{1}{\lambda} \left(1-\exp_g(-\lambda;\alpha,M)\right)
				\leqslant \frac{1}{\lambda}\left(1+\exp_g(-\lambda;\alpha,M)\right).
			\end{aligned}
		\end{displaymath}
	Furthermore,
		\begin{displaymath}
			\begin{aligned}
				\lim_{M\to \infty} \exp_g(-\lambda;\alpha,M)=&
				\lim_{M \to \infty} \exp\left(
				-\lambda \mu_g([\alpha,M)\setminus D_g)+\sum_{s\in [\alpha,M)\cap D_g} \log\left(
				1-\lambda\, \Delta g(s)
				\right)
				\right)\\
				=&\exp\left(\lim_{M\to \infty}\left[
				-\lambda \mu_g([\alpha,M)\setminus D_g)+\sum_{s\in [\alpha,M)\cap D_g} \log\left(
				1-\lambda\, \Delta g(s)
				\right)\right]
				\right)\\
				=& C(\lambda)\in [0,1].
			\end{aligned}
		\end{displaymath}
	Hence,
		\begin{displaymath}
			\left\lVert \exp_g(-\lambda;\alpha,\cdot)\right\rVert^p_{L_g^p([\alpha,\infty))} \leqslant 
			\frac{1}{\lambda}\left(1+C(\lambda)\right),\end{displaymath} and, therefore, \begin{displaymath} \left\lVert \exp_g(-\lambda;\alpha,\cdot)\right\rVert_{L_g^p([\alpha,\infty))} 
			\leqslant C(\lambda) \in [1,\infty).\qedhere
		\end{displaymath}
\end{proof}

\begin{rem} \label{lambdaelec} In Lemma~\ref{gexpgext1} it is crucial to take $\lambda >0$ such that
	$1-\lambda\, \Delta g(t)\in (0,1)$, $\forall t\in [\alpha,\infty)\cap D_g$. We observe that this choice is always possible. Indeed, we know that the set $\{s\in [\alpha,\infty)\cap D_g:\; \Delta g(s)\geq {1}/{2}\}$ is finite. If it is empty, we define $\lambda=1$, otherwise, we consider
	\begin{displaymath}
		\lambda=\frac{1}{2\, \max\{\Delta g(t):\; t\in [\alpha,\infty)\cap D_g, \; \Delta g(t)\geq {1}/{2}\}},
	\end{displaymath}
we have that, for any $t\in[\alpha,\infty)$,
either $t\notin \{s\in [\alpha,\infty)\cap D_g:\; \Delta g(s)\geq {1}/{2}\}$, in which case
		\begin{displaymath}
			\lambda \, \Delta g(t)< \frac{1}{4\, \max\{\Delta g(s):\; s\in [\alpha,\infty)\cap D_g, \; \Delta g(s)\geq {1}/{2}\}}\leqslant \frac{1}{2},
		\end{displaymath}
or $t\in \{s\in [\alpha,\infty)\cap D_g:\; \Delta g(s)\geq {1}/{2}\}$, and
		\begin{displaymath}
			\lambda \, \Delta g(t)\leqslant \lambda\, \max\{\Delta g(s):\; s\in [\alpha,\infty)\cap D_g, \; \Delta g(s)\geq {1}/{2}\}=\frac{1}{2}.
		\end{displaymath}
In any case, $\lambda \, \Delta g(t)\in (0,1/2)$, $\forall t\in [\alpha,\infty)\cap D_g$, and, in particular, 
	$1-\lambda\, \Delta g(t)\in (1/2,1)$, for every $ t\in [\alpha,\infty)\cap D_g$.
\end{rem}

We now prove the analogous result for the interval $(-\infty,\alpha)$.

\begin{lem} \label{gexpgext2} Let $\lambda>0$. We have that $\exp_g(\lambda;\alpha,\cdot)\in 
	W_g^{1,p}((-\infty,\alpha))$, for every $p\in [1,\infty]$. Furthermore, 
	$\left\lVert u\right\rVert_{W_g^{1,p}((-\infty,\alpha))}\leqslant C(\lambda)$, where $C(\lambda)>0$ is a constant that does not depend on $p$.
\end{lem}

\begin{proof} On the one hand,
	\begin{displaymath}
		\begin{aligned}
			\exp_g(\lambda;\alpha,t)=&\exp \left(
			-\lambda \mu_g([t,\alpha)\setminus D_g) + 
			\sum_{s\in [t,\alpha)\cap D_g} \log \left(
			1-\frac{\lambda\, \Delta g(s)}{1+\lambda\, \Delta g(s)}
			\right)
			\right)\\
			=&\exp \left(
			-\lambda \mu_g([t,\alpha)\setminus D_g) + 
			\sum_{s\in [t,\alpha)\cap D_g} \log \left(
			\frac{1}{1+\lambda\, \Delta g(s)}
			\right)
			\right)\\
			=&\exp \left(
			-\lambda \mu_g([t,\alpha)\setminus D_g) - 
			\sum_{s\in [t,\alpha)\cap D_g} \log \left(1+\lambda\, \Delta g(s)
			\right)
			\right).
		\end{aligned}
	\end{displaymath}
	Hence, it is clear that $\exp_g(\lambda;\alpha,t) \in 
	(0,1],\; \forall t\in (-\infty,\alpha)$. 

Now, for $x<y<\alpha$,
	\begin{displaymath}
		\exp_g(\lambda;\alpha,y)-\exp_g(\lambda;\alpha,x)=
		\int_x^y \lambda \, \exp_g(\lambda;\alpha,s)\, \operatorname{d}\mu_g(s).
	\end{displaymath}
	Thus, as in the previous case, it is enough to show that $\exp_g(\lambda;\alpha,\cdot)\in L_g^p((-\infty,\alpha))$, for every $p\in [1,\infty)$. We have again two cases. The case $p=\infty$ is obvious given that $\exp_g(\lambda;\alpha,t)\in (0,1]$, $\forall t\in (-\infty,\alpha)$, whereas for $1\leqslant p <\infty$,
\begin{displaymath}
\int_{(-\infty,\alpha)} \left\lvert \exp_g(\lambda;\alpha,s)\right\rvert^p\; \operatorname{d}\mu_g(s)=
\lim_{M \to -\infty} \int_{M}^{\alpha} \left\lvert \exp_g(\lambda;\alpha,s)\right\rvert^p\, \operatorname{d}\mu_g(s).
\end{displaymath}
Now, since
\begin{displaymath}
\exp_g(\lambda;\alpha,t)^p \leqslant \exp_g(\lambda;\alpha,t),\; \forall t\in (-\infty,\alpha),\; \forall 
p\in [1,\infty),
\end{displaymath}
we have that
\begin{displaymath}
\begin{aligned}
\int_{M}^{\alpha} \exp_g(\lambda;\alpha,s)^p\, \operatorname{d}\mu_g(s) \leqslant &
\int_{M}^{\alpha} \exp_g(\lambda;\alpha,s)\, \operatorname{d}\mu_g(s)
=\frac{1}{\lambda} \int_{M}^{\alpha}\lambda\, \exp_g(\lambda;\alpha,s)\, \operatorname{d}\mu_g(s)\\
=&\frac{1}{\lambda} \left(1-\exp_g(\lambda;\alpha,M)\right)
\leqslant \frac{1}{\lambda} \left(1+\exp_g(\lambda;\alpha,M)\right).
\end{aligned}
\end{displaymath}
On the other hand,
\begin{displaymath}
\begin{aligned}
\lim_{M \to -\infty} \exp_g(\lambda;\alpha,M)=&\lim_{M \to -\infty} 
\exp \left(-\lambda \mu_g([M,\alpha)\setminus D_g)+ 
\sum_{s\in [M,\alpha)\cap D_g} 
\log \left(1-\frac{\lambda\, \Delta g(t)}{1+\lambda\, \Delta g(t)}\right) \right)\\
=&\exp\left( 
\lim_{M \to -\infty} \left[ 
-\lambda \mu_g([M,\alpha)\setminus D_g)- 
\sum_{s\in [M,\alpha)\cap D_g} \log(1+\lambda \Delta g(t))
\right]
\right)\\
=& C(\lambda) \in [0,1].
\end{aligned}
\end{displaymath}
Hence,
		\begin{displaymath}
			\left\lVert \exp_g(\lambda;\alpha,\cdot)\right\rVert^p_{L_g^p((-\infty,\alpha))} \leqslant \frac{1}{\lambda}\left(1+C(\lambda) \right), \end{displaymath}
		and, therefore, \begin{displaymath} \left\lVert \exp_g(\lambda;\alpha,\cdot)\right\rVert_{L_g^p((-\infty,\alpha))} \leqslant C(\lambda)\in [1,\infty).\qedhere
		\end{displaymath}
The case $x\leqslant\alpha\leqslant y$ can be obtained by combining the previous two cases.
\end{proof}

In the following lemma we will prove that it is possible to concatenate 
two functions in $W_g^{1,p}$ under certain conditions. This will allow us to construct an extension operator.

\begin{lem} \label{concat} Let $u_1\in W_g^{1,p}(I_1)$ y $u_2\in W_g^{1,p}(I_2)$ be such that
	$\overline{I}_1\cap\overline{I}_2=\{\alpha\} \subset \mathbb{R}$ and $x<y$, $\forall x\in I_1$, 
	$\forall y\in I_2$ (observe that
	$I=I_1\cup I_2$ is an interval open from the right and closed from the left). Assume that 
 $u_1(\alpha)=u_2(\alpha)$ and define
	\begin{displaymath}
		u:t\in \overline{I} \rightarrow u(t)=\left\{
		\begin{aligned}
			& u_1(t),\; t\in I_1, \\
			& u_2(t),\; t\in I_2,
		\end{aligned}
		\right.\qquad
		\widetilde{u}:t\in \overline{I} \rightarrow u(t)=\left\{
		\begin{aligned}
			& \widetilde{u}_1(t),\; t\in I_1, \\
			& \widetilde{u}_2(t),\; t\in I_2,
		\end{aligned}
		\right.
	\end{displaymath}
where $\widetilde u_1$ and $\widetilde u_2$ are given by the definition of $W_g^{1,p}(I_1)$ and $W_g^{1,p}(I_2)$ respectively. Then, $u,\widetilde{u}\in L_g^p(I)$ and
	\begin{displaymath}
		u(y)-u(x)=\int_x^y \widetilde{u}(s)\, \operatorname{d}\mu_g(s),\; \forall x,y\in \overline{I}.
	\end{displaymath}
In particular, $u\in W_g^{1,p}(I)$ and
	\begin{displaymath}
		\left\lVert u\right\rVert_{W_g^{1,p}(I)}\leqslant \left\lVert u_1\right\rVert_{W_g^{1,p}(I_1)}+\left\lVert u_2\right\rVert_{W_g^{1,p}(I_2)}.
	\end{displaymath}
\end{lem}

\begin{proof} On the one hand,
	\begin{displaymath}
		\left\lVert u\right\rVert_{L_g^p(I)}^p = \int_I \left\lvert u(s)\right\rvert^p\, \operatorname{d}\mu_g(s)=\int_{I_1} \left\lvert u_1(s)\right\rvert^p\, \operatorname{d}\mu_g(s)+
		\int_{I_2} \left\lvert u_2(s)\right\rvert^p\, \operatorname{d}\mu_g(s),
	\end{displaymath}
	whence
	\begin{displaymath}
		\left\lVert u\right\rVert_{L_g^p(I)} = \left(\left\lVert u_1\right\rVert_{L_g^p(I_1)}^p+\left\lVert u_2\right\rVert_{L_g^p(I_2)}^p \right)^{1/p} 
		\leqslant \left\lVert u_1\right\rVert_{L_g^p(I_1)}+\left\lVert u_2\right\rVert_{L_g^p(I_2)}.
	\end{displaymath}
The case of $\widetilde{u}$ is analogous.

	On the other hand, given $x,y\in \overline{I}$, we have that, whether $x,y\in \overline{I}_1$ or $x,y\in \overline{I}_2$, 
it holds that
	\begin{displaymath}
		u(y)-u(x)=\int_x^y \widetilde{u}(s)\, \operatorname{d}\mu_g(s).
	\end{displaymath}
	Let us study the case $x\in I_1$ and $y\in I_2$. Then,
	\begin{align*}
			u(y)-u(x)=&u(y)-u(\alpha)+u(\alpha)-u(x)\\
			=&u_2(y)-u_2(\alpha)+u_1(\alpha)-u_1(x)\\
			=&\int_x^{\alpha}\widetilde{u}_1(s)\, \operatorname{d}\mu_g(s) +
			\int_{\alpha}^y \widetilde{u}_2(s)\, \operatorname{d}\mu_g(s) 
			=\int_x^y \widetilde{u}(s)\, \operatorname{d}\mu_g(s).\qedhere
	\end{align*}
\end{proof}

Now, we will use the previous results to show that, for
$a<b$, it is possible to extend functions in $W_g^{1,p}([a,b))$ to functions in $W_g^{1,p}(\mathbb{R})$. 
In particular, we have the following result.

\begin{thm}[Extension operator] \label{extth} Let $a<b$ and $1\leqslant p \leqslant \infty$. Then there exists a continuous linear operator $P:W_g^{1,p}([a,b)) \rightarrow W_g^{1,p}(\mathbb{R})$, 
	called the \emph{extension operator}, such that
	\begin{enumerate}
		\item $P f|_{[a,b]}=f$, $\forall f \in W_g^{1, p}([a,b))$,
		\item $\left\lVert P f\right\rVert_{L^p_g(\mathbb{R})} \leqslant \widetilde{C}\left\lVert f\right\rVert_{L^p_g([a,b))}$, $\forall f \in W^{1, p}_g([a,b))$,
		\item $\left\lVert P f\right\rVert_{W^{1, p}_g(\mathbb{R})} \leqslant \widetilde{C}\left\lVert f\right\rVert_{W_g^{1, p}([a,b))}$, $\forall 
		f \in W^{1, p}_g([a,b))$,
		\item $\left\lVert Pf\right\rVert_{0} \leqslant \widetilde{C} 
		\left\lVert f\right\rVert_{0}$, $\forall 
		f \in W^{1, p}_g([a,b))$.
	\end{enumerate}
	where $\widetilde{C}>0$ is a constant depending only on
	$p$ and $\mu_g([a,b))$. 
\end{thm}

\begin{proof} Given $f\in W_g^{1,p}([a,b))$, we define
	\begin{displaymath}
		\begin{aligned}
			P^-f:x\in (-\infty,a] & \longrightarrow & P^-f(x)=f(a)\, \exp_g(+\lambda^-;a,t)\in {\mathbb F} ,\\
			P^+f:x\in [b,+\infty) & \longrightarrow & P^+f(x)=f(b)\, \exp_g(-\lambda^+;b,t)\in {\mathbb F} ,
		\end{aligned}
	\end{displaymath}
	where $\lambda^->0$ is any element of choice, for instance, we can take $\lambda^-=1$, and $\lambda^+>0$ given by Remark~\ref{lambdaelec}. 
	Thanks to Lemmas ~\ref{gexpgext1} and~\ref{gexpgext2}, we have that 
	$P^-f\in W_g^{1,p}((-\infty,a))$ and $P^+f\in W_g^{1,p}([b,\infty))$. Furthermore, 
	it is clear that $P^-f(a)=f(a)$, $P^+f(b)=f(b)$ and
	\begin{displaymath}
		\begin{aligned}
			\left\lVert P^-f\right\rVert_{W_g^{1,p}((-\infty,a))}\leqslant & \left\lvert f(a)\right\rvert\, C(\lambda^-), \\
			\left\lVert P^+f\right\rVert_{W_g^{1,p}([b,\infty))}\leqslant & \left\lvert f(b)\right\rvert\, C(\lambda^+).
		\end{aligned}
	\end{displaymath} 

Let us consider $P:f\in W^{1,p}_g([a,b))\rightarrow Pf$ given by
	\begin{displaymath}
		Pf(x)=\begin{dcases}
		 P^-f(x),& x\in (-\infty,a), \\
			 f(x), & x\in [a,b), \\
			 P^+f(x), & x\in [b,\infty).
		\end{dcases}
	\end{displaymath}

By Lemma~\ref{concat}, we have that $Pf\in W_g^{1,p}(\mathbb{R})$. By this property and the definition of $P$, it is clear that property 1 in the statement of the theorem holds. Furthermore,
	\begin{displaymath}
		\begin{aligned}
			\left\lVert Pf\right\rVert_{W_g^{1,p}(\mathbb{R})}\leqslant & \left\lVert f\right\rVert_{W_g^{1,p}([a,b))}+ \left\lvert f(a)\right\rvert\, C(\lambda^-)+ \left\lvert f(b)\right\rvert\, C(\lambda^+)\\
			\leqslant & \left\lVert f\right\rVert_{W_g^{1,p}([a,b))}\left(1+C(p,\mu_g([a,b)))\,C(\lambda^-)+C(p,\mu_g([a,b)))\,C(\lambda^+)\right)\\
			=&\widetilde{C}(p,\mu_g([a,b)))\, \left\lVert f\right\rVert_{W_g^{1,p}([a,b))},
		\end{aligned}
	\end{displaymath}
where
	\begin{displaymath}
		\widetilde{C}(p,\mu_g([a,b)))=\left(1+C(p,\mu_g([a,b)))\,C(\lambda^-)+C(p,\mu_g([a,b)))\,C(\lambda^+)\right),
	\end{displaymath}
	and $C(p,\mu_g([a,b)))$ is the embedding constant of
	de $W^{1,p}_g([a,b))$ into $\operatorname{BC}_g([a,b])$. Therefore, we have proven properties 2 and 3.

Finally, given $t\in \mathbb{R}$, $\left\lvert Pf(t)\right\rvert=\left\lvert f(t)\right\rvert\leqslant \left\lVert f\right\rVert_{0}$, 
	if $t\in [a,b]$ and, in the case $t\notin [a,b]$, $\left\lvert Pf(t)\right\rvert\leqslant \max\{\left\lvert f(a)\right\rvert,\left\lvert f(b)\right\rvert\}\leqslant \left\lVert f\right\rVert_{0}$, so property 4 holds as well. 
\end{proof}

Let us see now that $W_g^{1,1}(I)$ is compactly embedded into 
$L^p_g(I)$ for every $1\leqslant p <\infty$.

\begin{thm}[Compactness of $W_g^{1,1}([a,b))$ into $L^p_g([a,b))$] 
The embedding of $W_g^{1,1}([a,b))$ into $L^p_g([a,b))$ is compact.
\end{thm}
\begin{proof}
Fix $p\in [1,\infty)$.	Let us consider $a,b\in \mathbb{R}$ such that $a<b$ and define
	\begin{displaymath}
		\mathcal{H}:=\{f \in W_g^{1,p}([a,b)):\; \left\lVert f\right\rVert_{W^{1,p}_g(a,b)}\leqslant 1\}.
	\end{displaymath}
 To prove the thesis of the theorem, it is enough to show that $\mathcal{H}$ is a relatively compact subset of 	$L^p_g([a,b))$.

	 For convenience, we may assume $g_c(\mathbb{R})=\mathbb{R}$ and $\sum_{t\in D_g}\Delta g(t)<\infty$. Otherwise, we can always construct $\widetilde{g}:\mathbb{R}\rightarrow 
	\mathbb{R}$ such that $\widetilde{g}(t)=g(t)$, $\forall t\in [a,b]$, $\widetilde{g}_c(\mathbb{R})=
	\mathbb{R}$ and $\sum_{t\in D_{\widetilde{g}}}\Delta \widetilde{g}(t)<\infty$.

	Let $I=[a,b)$, $P$ the extension operator given by Theorem~\ref{extth} and $\mathcal{G}=P(\mathcal{H})$. Let us show that $\mathcal{G}$ is a relatively compact subset of $L_g^p(\mathbb{R})$. In order to achieve this goal, it is enough to check that the conditions of Theorem~\ref{k-r-g} hold. In order to simplify the notation, we will denote by $f$ both the function defined in $I$ 
itself and its extension as given by Theorem~\ref{extth}.

Condition 1 in Theorem~\ref{k-r-g} states that $\{f(d_n):f\in\mathcal{G}\}\subset\mathbb{R}$ is bounded for each $d_n\in D_g$. 
Indeed, by Theorem~\ref{extth}, point 4, $\left\lvert f(d_n)\right\rvert\leqslant \left\lVert f\right\rVert_{0}\leqslant C\, \left\lVert f\right\rVert_{W^{1,1}_g(I)}\leqslant C$, for every $n\in \mathbb{N}$ and $f\in \mathcal{G}$ where $C$ is the embedding constant of $W^{1,1}_g([a,b))$ into $\operatorname{BC}_g([a,b])$.


Condition 2 in 	Theorem~\ref{k-r-g} can be checked taking into account that $\left\lvert f(d_n)\right\rvert\leqslant C$, for every $n\in \mathbb{N}$ and $f\in \mathcal{G}$, and $\sum_{t\in D_g} \Delta g(t)<\infty$. Indeed,
\begin{displaymath}
\sum_{k=1}^{\infty}\left\lvert f(d_k)\right\rvert^{p} \Delta g(d_k)\leqslant C^p \sum_{k=1}^{\infty} \Delta g(d_k)<\infty.
\end{displaymath}
Therefore, given $\varepsilon>0$, there exists $n\in \mathbb{N}$ such that
\begin{displaymath}
\sum_{k>n} \left\lvert f(d_k)\right\rvert^{p} \Delta g(d_k) < \varepsilon^p, \; \forall f\in\mathcal{G}.
\end{displaymath}

 Condition 3 in Theorem~\ref{k-r-g} can be checked by observing that, given $R>0$ and $f\in \mathcal{G}$
\begin{displaymath}
			\int_{A_R} \left\lvert f\right\rvert^p \, \operatorname{d}\mu_g= \sum_{t\in D_g\cap A_R} \left\lvert f(t)\right\rvert^p\, \Delta g(t)+
			\int_{g_c(A_R)} \left\lvert f\circ \gamma\right\rvert^p\, \operatorname{d} x,
		\end{displaymath}
		where $A_R=\{x\in \mathbb{R}:\; \left\lvert x\right\rvert>R\}$ and $\gamma$ is the pseudo-inverse of $g_c$. Now, 
	given $R$ sufficiently large,
\begin{displaymath}
\begin{aligned}
\int_{A_R} \left\lvert f\right\rvert^p\, \operatorname{d}\mu_g = & 
\lim_{M \to \infty} \int_{R}^{M} \left\lvert f(b)\right\rvert^p \, \exp_g(-\lambda^+; b,t)^p\, \operatorname{d}\mu_g(t)\\
& + \lim_{M \to -\infty} \int_{M}^{-R} \left\lvert f(a)\right\rvert^p\, \exp_g(+\lambda^-;a,t)^p\, \operatorname{d}\mu_g(t)\\
\leqslant & C^p\, \bigg( \lim_{M \to \infty} \int_{R}^{M} \exp_g(-\lambda^+; b,t)
\, \operatorname{d}\mu_g(t) \\ 
& +\lim_{M \to -\infty} \int_{M}^{-R} \exp_g(+\lambda^-; a,t)
\, \operatorname{d}\mu_g(t) \bigg)\\
=&C^p\, \bigg( \frac{1}{\lambda^+}\,\lim_{M \to \infty} \left[
\exp_g(-\lambda^+; b,R)-\exp_g(-\lambda^+; b,M)
\right]\\
&+\frac{1}{\lambda^-}\, \lim_{M \to -\infty} \left[
\exp_g(+\lambda^-; a,-R)- \exp_g(+\lambda^-; a,M)
\right]\bigg).
\end{aligned}
\end{displaymath}
		In particular, taking into account the hypotheses on $g_c$, 
		\begin{displaymath}
			\lim_{M\to \infty} \exp_g(-\lambda^+; b,M)=
			\lim_{M\to -\infty} \exp_g(+\lambda^-; a,t)=0,
		\end{displaymath}
	from which it follows that for every $\varepsilon>0$ there exists $R>0$ such that
		\begin{equation}\label{ccon1}
				\int_{A_R} \left\lvert f\right\rvert^p\, \operatorname{d}\mu_g \leqslant \varepsilon^p,\; \forall f\in \mathcal{G}.
		\end{equation}
		Now, $\{y\in \mathbb{R}:\; \left\lvert y\right\rvert> \max\{\left\lvert g_c(R)\right\rvert,\left\lvert g_c(-R)\right\rvert\}\}\subset\{g_c(x):\; \left\lvert x\right\rvert>R\}$, 
		whence, writing $\widetilde{R}:=\max\{\left\lvert g_c(R)\right\rvert,\left\lvert g_c(-R)\right\rvert\}$, we derive that
		\begin{displaymath}
			\int_{\left\lvert x\right\rvert>\widetilde{R}} \left\lvert f\circ \gamma(x)\right\rvert^p\, \operatorname{d} x \leqslant 
			\int_{g_c(A_R)} \left\lvert f\circ \gamma(x)\right\rvert^p\, \operatorname{d} x \leqslant \int_{A_R} \left\lvert f\right\rvert^p\, \operatorname{d}\mu_g(t),
		\end{displaymath}
	so 3 follows from~\eqref{ccon1}.

	Finally, we show that	Theorem~\ref{k-r-g}, point 4, holds. Given $f\in \mathcal{G}$ and $h>0$ (the case $h<0$ is analogous),
\begin{displaymath}
\begin{aligned}
f\circ\gamma(x+h)-f\circ\gamma(x)
=& \int_{\gamma(x)}^{\gamma(x+h)} f_g'(t)\, \operatorname{d}\mu_g(t)\\
=& \int_{\gamma(x)}^{\gamma(x+h)} f_g'(t)\, \operatorname{d}\mu_{g_C}(t)
+\int_{\gamma(x)}^{\gamma(x+h)} f_g'(t)\, \operatorname{d}\mu_{g_B}(t).
\end{aligned}
\end{displaymath}
Now, on the one hand,
\begin{displaymath}
\begin{aligned}
\int_{\gamma(x)}^{\gamma(x+h)} f_g'(t)\, \operatorname{d}\mu_{g_C}(t)
=& \int_{g_c([\gamma(x),\gamma(x+h)))} f_g'\circ \gamma (t)\, \operatorname{d} t \\
=& \int_x^{x+h} f_g'\circ \gamma (t)\, \operatorname{d} t\\ 
=& \,h \, \int_0^1 f_g'\circ \gamma (x+sh)\, \operatorname{d} s,
\end{aligned}
\end{displaymath}
so
\begin{displaymath}
			\left\lvert f\circ\gamma(x+h)-f\circ\gamma(x)\right\rvert\leqslant h \, \int_0^1 \left\lvert  f_g'\circ \gamma (x+sh)\right\rvert\, \operatorname{d} s
			+\int_{[\gamma(x),\gamma(x+h))} \left\lvert f_g'(t)\right\rvert\, \operatorname{d}\mu_{g_B}(t).
		\end{displaymath}
		If we integrate the first term on $\mathbb{R}$, 
		\begin{displaymath}
			\int_{\mathbb{R}}\bigg(
			\int_0^1 \left\lvert  f_g'\circ \gamma (x+sh)\right\rvert\, \operatorname{d} s\bigg)\, \operatorname{d} x =\int_0^1\bigg(\int_{\mathbb{R}}
			\left\lvert  f_g'\circ \gamma (x+sh)\right\rvert\,\operatorname{d} x \bigg)\, \operatorname{d} s
			\leqslant \left\lVert f_g'\right\rVert_{L_g^1(\mathbb{R})}\leqslant C.
		\end{displaymath}
	For the second term we have that
		\begin{displaymath}
			\begin{aligned}
				\int_{\mathbb{R}} \bigg(\int_{[\gamma(x),\gamma(x+h))} \left\lvert f_g'(t)\right\rvert\, \operatorname{d}\mu_{g_B}(t)\bigg)\, \operatorname{d} x=&
				\int_{\mathbb{R}} \bigg( \int_{\mathbb{R}} \chi_{[\gamma(x),\gamma(x+h))}(t)\, 
				\left\lvert f_g'(t)\right\rvert\, \operatorname{d}\mu_{g_B}(t)\bigg)\, \operatorname{d} x.
			\end{aligned}
		\end{displaymath}
	Now,
		\begin{displaymath}
			t\in [\gamma(x),\gamma(x+h))\Leftrightarrow x\leqslant g_c(t)\leqslant x+h \Leftrightarrow g_c(t)-h\leqslant x\leqslant g_c(t),
		\end{displaymath}
	so
		\begin{displaymath}
			\begin{aligned}
				\int_{\mathbb{R}} \bigg( \int_{\mathbb{R}} \chi_{[\gamma(x),\gamma(x+h))}(t)\, 
				\left\lvert f_g'(t)\right\rvert\, \operatorname{d}\mu_{g_B}(t)\bigg)\, \operatorname{d} x =& \int_{\mathbb{R}} \bigg( \int_{\mathbb{R}} 
				\chi_{[g_c(t)-h,g_c(t)]}(x)\left\lvert f_g'(t)\right\rvert\, \operatorname{d}\mu_{g_B}(t)\bigg)\, \operatorname{d} x\\
				=& \int_{\mathbb{R}} \left\lvert f_g'(t)\right\rvert\, \bigg(\int_{\mathbb{R}} \chi_{[g_c(t)-h,g_c(t)]}(x)\, \operatorname{d} x \bigg)\, \operatorname{d}\mu_{g_B}(t)\\
				=& h \, \int_{\mathbb{R}} \left\lvert f_g'(t)\right\rvert\,\operatorname{d}\mu_{g_B}(t) 
				\leqslant h \, \left\lVert f_g'\right\rVert_{L_g^1(\mathbb{R})}\leqslant C\, h.
			\end{aligned}
		\end{displaymath}
	To summarize,
		\begin{displaymath}
			\int_{\mathbb{R}} \left\lvert f\circ\gamma(x+h)-f\circ\gamma(x)\right\rvert \operatorname{d} x\leqslant 2 C\, h.
		\end{displaymath}
	Finally,
		\begin{displaymath}
			\begin{aligned}
				\int_{\mathbb{R}} \left\lvert f\circ\gamma(x+h)-f\circ\gamma(x)\right\rvert^p \operatorname{d} x = &
				\int_{\mathbb{R}} \left\lvert f\circ\gamma(x+h)-f\circ\gamma(x)\right\rvert^{p-1} \, \left\lvert f\circ\gamma(x+h)-f\circ\gamma(x)\right\rvert \operatorname{d} x\\
				\leqslant & \left\lVert f\right\rVert_{0}^{p-1}\, \int_{\mathbb{R}} \left\lvert f\circ\gamma(x+h)-f\circ\gamma(x)\right\rvert \operatorname{d} x\leqslant 2C\, h,
			\end{aligned}
		\end{displaymath}
	where $C$ does not depend on $f$. Taking $\rho$ sufficiently small and $\left\lvert h\right\rvert<\rho$ we show that	Theorem~\ref{k-r-g}, point 4, holds.
\end{proof}

\section{An application to decomposable functions}

In this section we will apply the compactness results obtained in the previous sections to the space of those functions $f:[a,b]\to{\mathbb F}$ that can be expressed as the sum of a continuous function and a jump function. We characterize these functions in the following result.
\begin{lem}\label{lemeqbv}Let $f\in\operatorname{BC}_g([a,b],{\mathbb F})$ be such that $f(t^+)$ exists for every $t\in D_g\cap [a,b)$. The following are equivalent:
	\begin{enumerate}
		\item $\sum_{t\in D_g\cap[a,b)}\left\lvert \Delta f(t)\right\rvert<\infty$.
		\item $D_g\cap[a,b)$ is finite or, for every $\{t_n\}_{n\in{\mathbb N}}\subset D_g\cap[a,b)$, we have that $(f_n)_{n\in{\mathbb N}}$, where \[ f_n(t):=\sum\limits_{\substack{t_k<t\\  k=1,\dots,n}}\Delta f(t_k),\]  for $t\in[a,b)$ is a convergent sequence in $\operatorname{BC}_g([a,b],{\mathbb F})$.
		\item The function $h(t):=f'_g(t)$ if $t\in D_g\cap[a,b)$, $h(t)=0$ if $t\in[a,b]\backslash D_g$ is well defined and belongs to $L_g^1([a,b],{\mathbb F})$.
		\item The function $f^B(t)= \sum_{s\in[a,t)} \Delta f(s)$ is well defined, $f^B\in\operatorname{AC}_{g^B}([a,b],{\mathbb F})\cap\operatorname{AC}_g([a,b],{\mathbb F})$ and $f^C:=f-f^B\in \operatorname{BC}_g([a,b],{\mathbb F})\cap \operatorname{BC}([a,b],{\mathbb F})$.
	\end{enumerate}
\end{lem}
\begin{proof}
	Observe that we can talk of $\Delta f$ at the points of $D_g$ because, by hypothesis, $f(t^+)$ exists at the points of $D_g$. The case where $D_g\cap[a,b)$ is finite is straightforward, so we will deal only with the case where $D_g\cap[a,b)$ is infinite and countable.

	1$\Rightarrow$2. Let $\{t_n\}_{n\in{\mathbb N}}\subset D_g$, $\varepsilon\in{\mathbb R}^+$. There exists $N\in{\mathbb N}$	such that $\sum_{n>N}\left\lvert \Delta f(t_n)\right\rvert<\varepsilon$. Thus, for $m\geqslant n\geqslant N$, \[ \left\lVert f_m-f_n\right\rVert=\sup_{t\in[a,b]}\left\lvert \sum\limits_{\substack{t_k<t\\k=n+1,\dots,m}}\Delta f(t_k)\right\rvert\leqslant\sup_{X\subset\{n+1,\dots,m\}}\left\lvert \sum_{k\in X}\Delta f(t_k)\right\rvert\leqslant\sum_{n>N}\left\lvert \Delta f(t_n)\right\rvert<\varepsilon.\] 
	Hence, $(f_n)_{n\in{\mathbb N}}$ is a Cauchy sequence in $\operatorname{BC}_g([a,b],{\mathbb F})$ and, therefore, convergent.

	2$\Rightarrow$1. Let $D_g^+:=\{t\in D_g\ :\ \Delta f(t)>0\}=\{r_n\}_{n\in\Lambda_1}$, $D_g^-:=\{t\in D_g\ :\ \Delta f(t)<0\}=\{s_n\}_{n\in\Lambda_2}$, where $\Lambda_1,\Lambda_2\subset{\mathbb N}$. If $\Lambda_1$ is finite, it is clear that $\sum_{k\in\Lambda_1}\left\lvert \Delta f(r_k)\right\rvert<\infty$. Otherwise, the sequence $(f_n)_{n\in{\mathbb N}}$ where \[ f_n(t):=\sum\limits_{\substack{r_k<t\\k=1,\dots,n}}\Delta f(r_k)\]  is convergent in $\operatorname{BC}_g([a,b],{\mathbb F})$, so there exists \[ \lim_{n\to\infty}f_n(b)=\lim_{n\to\infty}\sum\limits_{\substack{r_k<b\\k=1,\dots,n}}\Delta f(r_k)=\lim_{n\to\infty}\sum\limits_{\substack{r_k<b\\k=1,\dots,n}}\left\lvert \Delta f(r_k)\right\rvert=\sum\limits_{r_k<b}\left\lvert \Delta f(r_k)\right\rvert=\sum\limits_{k\in\Lambda_1}\left\lvert \Delta f(r_k)\right\rvert<\infty.\] 
	If $\Lambda_2$ is finite, it is clear that $\sum_{k\in\Lambda_1}-\left\lvert \Delta f(s_k)\right\rvert>-\infty$. Otherwise, $\Lambda_2$ is infinite, the sequence $(f_n)_{n\in{\mathbb N}}$ where \[ f_n(t):=\sum\limits_{\substack{s_k<t\\k=1,\dots,n}}\Delta f(s_k)\]  is convergent, so there exists
	\begin{align*}\lim_{n\to\infty}f_n(b)= & \lim_{n\to\infty}\sum\limits_{\substack{s_k<t\\k=1,\dots,n}}\Delta f(s_k)=\lim_{n\to\infty}\sum\limits_{\substack{s_k<t\\k=1,\dots,n}}-\left\lvert \Delta f(s_k)\right\rvert=-\sum\limits_{s_k<t}\left\lvert \Delta f(s_k)\right\rvert\\ = & -\sum\limits_{k\in\Lambda_2}\left\lvert \Delta f(s_k)\right\rvert >-\infty.\end{align*}
	Combining these two facts, $\sum_{t\in D_g}\left\lvert \Delta f(t)\right\rvert<\infty$.

	1$\Rightarrow$3. $f'_g$ is well defined on the points of $D_g$ because $f(t^+)$ exists at the points of $D_g$. $h$ is $\mu_g$-measurable because $h=0$ in $[a,b]\backslash D_g$ and $D_g$ is countable and $g$-measurable.
	In fact, $\mu_g(f^{-1}(E))=\sum_{t\in f^{-1}(E)}\Delta g(t)$. Now it is enough to observe that
	\begin{align*}\int_{[a,b)}\left\lvert h(t)\right\rvert\operatorname{d} \mu_g(t)= & \int_{[a,b)\cap D_g}\left\lvert f_g'(t)\right\rvert\operatorname{d} \mu_g(t)=\int_{[a,b)\cap D_g}\frac{\left\lvert \Delta f(t)\right\rvert}{\Delta g(t)}\operatorname{d} \mu_g(t)\\= & \sum_{t\in D_g\cap[a,b)}\frac{\left\lvert \Delta f(t)\right\rvert}{\Delta g(t)}\Delta g(t)=\sum_{t\in D_g\cap[a,b)}\left\lvert \Delta f(t)\right\rvert<\infty.\end{align*}
	Therefore, $h\in L_g^1([a,b],{\mathbb F})$.

	3$\Rightarrow$1. Since $h\in L_g^1([a,b],{\mathbb F})$, we have that $\infty>\int_{[a,b)}\left\lvert h(t)\right\rvert\operatorname{d} \mu_g(t)=\sum_{t\in D_g\cap[a,b)}\left\lvert \Delta f(t)\right\rvert$.

	3$\Rightarrow$4. By definition of $h$,
	\begin{align*}\int_{[a,t)}h(s)\operatorname{d} \mu_g(s)= & \int_{[a,t)\cap D_g}f_g'(s)\operatorname{d} \mu_g(s)=\int_{[a,t)\cap D_g}\frac{\Delta f(s)}{\Delta g(s)}\operatorname{d} \mu_g(s)= \sum_{s\in D_g\cap[a,t)}\frac{\Delta f(s)}{\Delta g(s)}\Delta g(t)\\ = & \sum_{s\in D_g\cap[a,t)}\Delta f(s)=\sum_{s\in[a,t)}\Delta f(s)=f^B(t),\end{align*}
	so $f^B$ is well defined and, by the Fundamental Theorem of Calculus, $f^B\in\operatorname{AC}_g([a,b],{\mathbb F})$, which implies that $f^B$ is $\mu_{g^B}$-measurable. Furthermore,
	\[ f^B(t)=\int_{[a,t)\cap D_g}f_g'(s)\operatorname{d} \mu_g(s)=\int_{[a,t)}f_g'(s)\operatorname{d} \mu_{g^B}(s).\] 
	Hence, by the Fundamental Theorem of Calculus, $f^B\in\operatorname{AC}_{g^B}([a,b],{\mathbb F})$. Also, since 3$\Rightarrow$4, we have that $\sum_{t\in D_g\cap[a,b)}\left\lvert \Delta f(t)\right\rvert<\infty$. Therefore, for every $t\in\ [a,b]$,
	\[ -\infty<-\sum_{s\in D_g\cap[a,b)}\left\lvert \Delta f(s)\right\rvert<\sum_{s\in D_g\cap[a,t)}\Delta f(s)=f^B(t)<\sum_{s\in D_g\cap[a,b)}\left\lvert \Delta f(s)\right\rvert<\infty,\] 
	and we conclude that $f^B$ is bounded.

	Since $f$ and $f^B$ are $g$-continuous and bounded, so is $f^C=f-f^B$. Finally, $f^C$ is continuous. Indeed, since $f^C$ is $g$-continuous, it is enough to see what happens at the points of $D_g$. Let $t\in D_g$. Then
	\begin{align*}f^C(t^+)-f^C(t)= & f(t^+)-f^B(t^+)-f(t)+f^B(t)\\ = &\Delta f(t)-\lim_{r\to t^+}\sum_{s\in D_g\cap[a,r)}\Delta f(s)+\sum_{s\in D_g\cap[a,t)}\Delta f(s) \\ = & \Delta f(t)-\lim_{r\to t^+}\sum_{s\in D_g\cap[t,r)}\Delta f(s)=-\lim_{r\to t^+}\sum_{s\in D_g\cap(t,r)}\Delta f(s)=0.\end{align*}
	Therefore, $f^C$ is continuous.

	4$\Rightarrow$1. Just observe that
	\[ \infty>\int_{[a,b)}\left\lvert (f^B)'_g\right\rvert(s)\operatorname{d} \mu_{g^B}(s)=\sum_{t\in D_g\cap[a,b)}\left\lvert \Delta f(t)\right\rvert.\qedhere\] 
\end{proof}
\begin{rem}\label{remord}
	The Riemann series theorem states that a convergent series is conditionally convergent if and only if the terms of the series can be rearranged in such a way that the sum of the new series is any fixed real number (or $\pm\infty$). As a consequence of this theorem, a convergent series is absolutely convergent if and only if the series of any subsequence of the sequence of terms is convergent. Points 1 and 2 in Lemma~\ref{lemeqbv} are reminiscent of this result, so it is only natural to wonder whether we can prove a similar, but less restrictive, lemma where we do not deal with the absolute convergence of $\sum_{t\in D_g}\left\lvert \Delta f(t)\right\rvert$ (remember that $D_g$ is countable, so this sum can be interpreted as a series in the traditional sense), but with some notion of conditional convergence of $\sum_{t\in D_g}\Delta f(t)$ or, even better, of $\sum_{s<t}\Delta f(s)$ for every $t\in[a,b]$ (see Lemma~\ref{lemeqbv}.4). This is not in general possible. The conditional convergence of series arises from the well order of the natural numbers and, although $D_g$ (or $D_g\cap[a, t)$ for some $t$) inherits the total order of the real numbers, it is not in general a well order, so any definition of $\sum_{s<t}\Delta f(s)$ would rely on a particular choice of a well order $\{t_n\}_{n\in{\mathbb N}}= D_g$ and, for a different order, the sum may yield different results.
\end{rem}

\begin{dfn}We will denote by $\operatorname{DC}_g([a,b],{\mathbb F})$ the set of functions satisfying all of the conditions in Lemma~\ref{lemeqbv} and consider in this space the norm
	\[ \left\lVert f\right\rVert_{\operatorname{DC}}:=\left\lVert f\right\rVert_\infty+\left\lVert f_g'\chi_{D_g}\right\rVert_1.\] 
\end{dfn}
\begin{rem}Observe that, for $f\in\operatorname{DC}_g([a,b],{\mathbb F})$,
	\[ \left\lVert f_g'\chi_{D_g}\right\rVert_1=\int_{[a,t)\cap D_g}\left\lvert f_g'(s)\right\rvert\operatorname{d} \mu_g(s)= \sum_{s\in D_g\cap[a,t)}\frac{\left\lvert \Delta f(s)\right\rvert}{\Delta g(s)}\Delta g(t)= \sum_{s\in D_g\cap[a,t)}\left\lvert \Delta f(s)\right\rvert<\infty,\] 
	so the norm $\left\lVert \cdot \right\rVert_{\operatorname{DC}}$ is well defined.
\end{rem}

Now we will study those cases where the decomposition can be taken as a product instead of a sum. This necessitates of the following result.
\begin{lem}\label{lemgb}If $f\in \operatorname{AC}_{g^B}([a,b],{\mathbb F})$ then $f(t)=f(a)+ \sum_{s\in[a,t)} \Delta f(s)$ for every $t\in [a,b]$.
\end{lem}
\begin{proof}
	By the fundamental Theorem of Calculus,
	\[ f(t)-f(a)=\int_{[a,t)}f_{g^B}'(s)\operatorname{d} \mu_g^B(s)=\sum_{s\in[a,t)}f_{g^B}'(s)\Delta g^B(s)=\sum_{s\in[a,t)}\Delta f(s).\qedhere\] 
\end{proof}

Let $\ln:{\mathbb R}^+\to{\mathbb R}$ be the real logarithm, $\arg_\sigma:{\mathbb C}\backslash\{0\}\to[\sigma-\pi,\sigma+\pi)$ be the $\sigma$-branch of the argument function (that is, $z=\left\lvert z\right\rvert e^{i\arg_\sigma z}$) and let $\log_\sigma:{\mathbb C}\backslash\{0\}\to {\mathbb R}$ be the $\sigma$-branch of the logarithm (that is, $\log_\sigma z=\ln\left\lvert z\right\rvert+i\arg_\sigma z$). We denote by $\log$ the principal branch of the complex logarithm, that is $\log_0$.

\begin{pro}
	Let $f\in\operatorname{BC}_g([a,b],{\mathbb F})$ be such that for every $t\in D_g$ there exists $f(t^+)$. Define $D_{g,f}=\{t\in D_g\ :\ \Delta f(t)\ne 0\}$. Assume that for every $t\in D_{g,f}$ there exists $\delta\in(0,b-t)$ such that $f(s)\ne0$ for every $s\in[t,t+\delta)$. If there exists $\alpha\in{\mathbb R}$ such that $\sum_{t\in D_{g,f}}\lvert\log_\alpha \left( \lim_{s\to t^+}f(t)\right) -\log_\alpha f(t)\rvert<\infty$ then
	we can write $f=\varphi\psi$ where $\varphi\in\operatorname{AC}_{g}([a,b],{\mathbb F})\cap \operatorname{AC}_{g^B}([a,b],{\mathbb F})$ , $\psi\in\operatorname{BC}_{g}([a,b],{\mathbb F})\cap \operatorname{BC}([a,b],{\mathbb F})$ and, for $t\in[a,b]$,
	\begin{align*}\varphi(t)= & \exp\left( \sum_{s\in D_{g,f}\cap[a,t)} \left[ \log_\alpha \left( \lim_{r\to s^+}f(s)\right) -\log_\alpha f(s)\right]\right) = \prod_{s\in D_{g,f}\cap[a,t)}\left( 1+\frac{\Delta f(t)}{f(t)}\right) \\ = & 1+\sum_{s\in D_{g,f}\cap[a,t)} \frac{\Delta f(s)}{f(s)}\Delta\varphi(s).\end{align*}
\end{pro}
\begin{proof}

	Define \[ \varphi(t)=\exp\left( \sum_{s\in D_{g,f}\cap[a,t)} \left[ \log_\alpha \left( \lim_{r\to s^+}f(s)\right) -\log_\alpha f(s)\right]\right) \]  for $t\in[a,b]$. Let us check that $\varphi$ is $g^B$-absolutely continuous. Let 
	\[ h(t)=\begin{dcases}\frac{ \log_\alpha \left( \lim_{r\to s^+}f(t)\right) -\log_\alpha f(t)}{\Delta g(t)}, & t\in D_{g,f},\\ 0, & t\in [a,b]\backslash D_{g,f}. \end{dcases}\] 
	The $\mu_g$ and $\mu_{g^B}$ mesurability of this function is argued as in the proof of Proposition~\ref{lemeqbv}. Furthermore, we have that
	$h\in{ L}^1_{g^B}([a,b],{\mathbb C})\cap { L}^1_{g}([a,b],{\mathbb C})$ and 
	\[ \varphi(t)=\exp\left( \int_{[a,t)}h(s)\operatorname{d}\mu_{g^B}(s)\right) =\exp\left( \int_{[a,t)}h(s)\operatorname{d}\mu_{g}(s)\right) ,\] 
	so $\varphi\in\operatorname{AC}_{g}([a,b],{\mathbb F})\cap \operatorname{AC}_{g^B}([a,b],{\mathbb F})$. 

	Besides,
	\begin{align*}\varphi(t)= &\exp\left( \sum_{s\in D_{g,f}\cap[a,t)} \left[ \log_\alpha \left( \lim_{r\to s^+}f(s)\right) -\log_\alpha f(s)\right]\right) \\= & \prod_{s\in D_{g,f}\cap[a,t)} \exp\left[ \log_\alpha \left( \lim_{r\to s^+}f(s)\right) -\log_\alpha f(s)\right]\\ = & \prod_{s\in D_{g,f}\cap[a,t)}\frac{f(s^+)}{f(s)}=\prod_{s\in D_{g,f}\cap[a,t)}\left( 1+\frac{\Delta f(s)}{f(s)}\right) .\end{align*}

	On the other hand, by Lemma~\ref{lemgb}, $\varphi(t)=\varphi(a)+\sum_{s\in[a,t)}\Delta \varphi(s)$.	Observe that $\varphi(a)=1$. Furthermore, 
	\begin{align*}\Delta \varphi(t)= & \lim_{r\to t^+}\varphi(r)-\varphi(t)=\lim_{r\to t^+}\left( \frac{\varphi(r)}{\varphi(t)}-1\right) \varphi(t) \\ = & \lim_{r\to t^+}\left( \prod_{s\in D_{g,f}\cap[t,r)}\left( 1+\frac{\Delta f(s)}{f(s)}\right) -1\right) \varphi(t)\\ = & \left( \left( 1+\frac{\Delta f(t)}{f(t)}\right) -1\right) \varphi(t)= \frac{\Delta f(t)}{f(t)}\varphi(t).\end{align*}
	Therefore,
	\[ \varphi(t)= 1+\sum_{s\in[a,t)} \frac{\Delta f(s)}{f(s)}\varphi(s).\]  
	Observe that 
	\[ \varphi(t)\geqslant\exp\left( -\sum_{s\in D_{g,f}} \left\lvert  \log_\alpha \left( \lim_{r\to s^+}f(s)\right) -\log_\alpha f(s)\right\rvert\right) =:M>0.\] 
	Define $\psi(t)=f(t)/\varphi(t)$. $f,\varphi\in \operatorname{BC}_{g}([a,b],{\mathbb F})$, so $\psi\in \operatorname{C}_{g}([a,b],{\mathbb F})$ by \cite[Lemma 2.14]{Fernandez2022a}. Furthermore, since $\varphi(t)\geqslant M$, $\psi$ is bounded. Finally, if $t\in [a,b]\backslash D_g$, $\psi$ is continuous at $t$, so it is left to check what happens for $t\in D_g$. In that case,
	\begin{align*}\psi(t^+)-\psi(t)= & \frac{f(t^+)}{\varphi(t^+)}-\frac{f(t)}{\varphi(t)}=\lim_{r\to t^+}\left( f(r)\prod_{s\in D_{g,f}\cap[a,r)}\frac{f(s)}{f(s^+)}\right) - f(t)\prod_{s\in D_{g,f}\cap[a,t)}\frac{f(s)}{f(s^+)}\\ = & \left[\prod_{s\in D_{g,f}\cap[a,t)}\frac{f(s)}{f(s^+)}\right]\left[\lim_{r\to t^+}\left( f(r)\prod_{s\in D_{g,f}\cap[t,r)}\frac{f(s)}{f(s^+)}\right) - f(t)\right].
	\end{align*}
	Since
	\[ \lim_{r\to t^+}\left( f(r)\prod_{s\in D_{g,f}\cap[t,r)}\frac{f(s)}{f(s^+)}\right) =f(t^+)\frac{f(t)}{f(t^+)}=f(t),\] 
	we conclude that $\psi(t^+)-\psi(t)=0$ and $\psi$ is continuous.
\end{proof}

Finally, we obtain a compactness result for $\operatorname{DC}_g([a,b],{\mathbb F})$. We consider $D_g=\{d_k\}_{k\in\Lambda}$, $\Lambda\subset{\mathbb N}$.

\begin{thm}
	$S\subset \operatorname{DC}_g([a,b],{\mathbb F})$ is totally bounded if and only if
	\begin{enumerate}
		\item			$S(t)$ is bounded for all $t\in [a,b]$,
		\item		$S$ is $g$-equicontinuous,
		\item		$S$ is $g$-stable,
		\item for every $\varepsilon>0$ there exists $n$ such that for any $f\in S$,
		\[ \sum_{k>n}\left\lvert f(d_k)\right\rvert<\varepsilon.\] 
	\end{enumerate}
\end{thm}
\begin{proof}
	Total boundedness of $S$ is equivalent to the total boundedness of both $\widetilde{S}=\{f:f\in S\}\subset\operatorname{BC}_g([a,b],\mathbb{F})$ and $\widehat{S}=\{f_g'\chi_{D_g}:f\in S\}\subset{L}^1([a,b],\mathbb{F})$. \par
	Since $\gamma$ (the pseudoinverse of $g^C$) is strictly increasing, $\gamma^{-1}(D_g)$ is at most countable, so conditions $3$ and $4$ in Theorem~\ref{k-r-g} are satisfied (the functions involved are zero except on sets of null Lebesgue measure). Condition $2$ translates as
	\[ \sum_{k>n_\varepsilon}\left\lvert f_g'(d_k)\right\rvert\Delta g(d_k)=\sum_{k>n_\varepsilon}\left\lvert f(d_k)\right\rvert<\varepsilon,\] 
	for every $\varepsilon>0$, $f\in S$ and some $n_\varepsilon$, so applying Theorems~\ref{thmbc} and~\ref{k-r-g} we obtain the result. 
\end{proof}

\bibliography{TR}
\bibliographystyle{spmpsciper}

\end{document}